\documentclass[11pt,reqno]{amsart}
\usepackage{amsthm,amssymb,amsbsy,epsfig,graphicx,color}
\usepackage{hyperref}

\newcommand{\nc}{\newcommand}
\nc{\be}{\begin{equation}} \nc{\ee}{\end{equation}}
\nc{\bea}{\begin{eqnarray}} \nc{\eea}{\end{eqnarray}}
\nc{\disp}{\displaystyle} \nc{\ade}{\mbox{$A$-$D$-$E$}}
\nc{\calN}{{\cal N}} \nc{\calC}{{\cal C}} \nc{\calM}{{\cal M}}
\nc{\calS}{{\cal S}} \nc{\phit}{\hat{\varphi}}
\nc{\chit}{\hat{\chi}} \nc{\hcalN}{\hat{\calN}}
\nc{\hcalS}{\hat{\calS}} \nc{\hS}{\hat{S}}
\nc{\sigmad}{\sigma^\dagger} \nc{\psid}{\psi^\dagger}

\definecolor{IndianRed}{rgb}{0.8,0.36,0.36}

\newtheorem{lemma}{Lemma}

\newtheorem{conj}{Conjecture}
\newtheorem{defn}{Definition}
\newtheorem{prop}{Proposition}

\newtheorem{example}{Example}

\font\tenmsb=msbm10\font\sevenmsb=msbm7 \font\fivemsb=msbm5
\newfam\msbfam
\textfont\msbfam=\tenmsb \scriptfont\msbfam=\sevenmsb
\scriptscriptfont\msbfam=\fivemsb

\def\binom#1#2{{#1\choose #2}}

\def\ket#1{|#1\rangle}
\def\e{{\rm e}}
\def\d{{\rm d}}
\def\i{{\rm i}}

\def\T{\mathcal{T}}

\def\D{\mathcal{D}}

\numberwithin{equation}{section}

\begin{document}
\begin{titlepage}
\title[Plane partitions, $q$KZ and Hirota]{Punctured plane partitions and the
$q$-deformed Knizhnik--Zamolodchikov and Hirota equations}

\author{Jan de Gier}
\address{Jan de Gier,
Department of Mathematics and Statistics, The University
of Melbourne, VIC 3010, Australia}
\email{degier@ms.unimelb.edu.au}
\thanks{JdG gratefully acknowledges support of the Australian Research Council.}

\author{Pavel Pyatov}
\address{Pavel Pyatov,
Max Planck Institute for Mathematics,
Vivatsgasse 7, D-53111 Bonn, Germany \&
Bogoliubov Laboratory of Theoretical Physics, Joint Institute
for Nuclear Research, 141980 Dubna, Moscow Region, Russia}
\email{pyatov@theor.jinr.ru}
\thanks{
PP was supported by the
DFG-RFBR grant 436 RUS 113/909/0- 1(R)
and 07-02-91561-a, and by the grant of the Heisenberg-Landau foundation.}

\author{Paul Zinn-Justin}
\address{Paul Zinn-Justin,
LPTMS (CNRS, UMR 8626), Univ Paris-Sud, 91405 Orsay, France; and
LPTHE (CNRS, UMR 7589), Univ Pierre et Marie Curie-Paris6 and Univ
Denis Diderot-Paris7, 75252 Paris, France.}
\email{pzinn@lptms.u-psud.fr}
\thanks{PZJ was supported by
EU Marie Curie Research Training Networks
``ENRAGE'' MRTN-CT-2004-005616, ``ENIGMA'' MRT-CT-2004-5652,
ESF program ``MISGAM''
and ANR program ``GIMP'' ANR-05-BLAN-0029-01.}

\date{\today}

\begin{abstract}
We consider partial sum rules for the homogeneous limit of the
solution of the $q$-deformed Knizhnik--Zamolodchikov
equation with reflecting boundaries
in the Dyck path representation. We show that these partial
sums arise in a solution of the discrete Hirota equation, and
prove that they are the generating functions of $\tau^2$-weighted
punctured cyclically symmetric transpose complement plane
partitions where $\tau=-(q+q^{-1})$. In the cases of
no or minimal punctures, we prove that these generating
functions coincide with $\tau^2$-enumerations of vertically
symmetric alternating sign matrices and modifications thereof.
\end{abstract}

\maketitle

{\footnotesize\tableofcontents}
\end{titlepage}

\section{Introduction}

The discrete Hirota equation \cite{Hirota,Zabrodin} is a
difference equation for a function $f$ in three discrete
variables:
\begin{multline}
f(n,i,j)f(n-2,i,j)=f(n-1,i-1,j)f(n-1,i+1,j)\\+\tau^2f(n-1,i,j-1)f(n-1,i,j+1).
\label{Hirota}
\end{multline}
The factor $\tau^2$ can be absorbed into $f$ by rescaling $f(n,i,j) \rightarrow
\tau^{-j(j+1)} f(n,i,j)$. Equation \eqref{Hirota} is also known as
the octahedron recurrence as the variables $n,i,j$ are natural
coordinates for the vertices of an octahedron, see e.g.
\cite{Speyer}.

In order to fix a particular solution to \eqref{Hirota} we need to
specify boundary conditions. An interesting choice
is the following. Let $A=(a_{ij})_{1\leq i,j\leq n}$ be
an $n\times n$ matrix and set
\be
f(0,i,j)=1,\qquad f(1,i,j)= a_{ij}\quad 1\leq i,j\leq n.
\label{ASMboundary}
\ee
With these boundary conditions, the discrete Hirota equation
defines the $\tau^2$-deformation of the determinant studied by
Robbins and Rumsey \cite{RobbinsR}. In particular, the
$\tau^2$-determinant of $A$ is defined by
\be
|A|_{\tau^2} = f(n,1,1).
\ee
For $\tau=\i$ the $\tau^2$-determinant reduces to the ordinary
determinant which can be expanded as a sum over permutation
matrices. In their famous work \cite{RobbinsR}, Robbins and Rumsey
showed that solutions to \eqref{Hirota} with the boundary
condition \eqref{ASMboundary} can be written as a sum over
alternating sign matrices (ASMs), see also \cite{Bressoud}.

It is well known that ASMs are equinumerous to totally symmetric
self-complementary plane partitions (TSSCPPs). Surprisingly,
generating functions of $\tau^2$-enumerations of TSSCPPs and other
symmetry classes of plane partitions, as studied by Robbins
\cite{Robbins}, appear as normalisations for homogeneous solutions
of the $q$-deformed Knizhnik--Zamolod\-chikov ($q$KZ) equation
recently obtained by Di Francesco and Zinn-Justin
\cite{DF07,DFZJ07}. This result will be generalised below, when we
consider a class of punctured cyclically symmetric transpose
complement plane partitions (PCSTCPPs), whose weighted
enumerations also arise in the $q$KZ equation.

It was already observed by Robbins and Kuperberg
\cite{Robbins,Kuperberg} that the $\tau^2$-enumerations of CSTCPPs
(without puncture) are closely related to $\tau^2$ enumerations of
vertically symmetric alternating-sign matrices (VSASMs) and other
symmetry classes of ASMs. A precise statement will be proved
below. Closing the circle, these enumerations comprise a
particular solution of \eqref{Hirota}, albeit with a different
boundary condition than \eqref{ASMboundary}. We hope that this
paper will be a further step in resolving the Razumov--Stroganov
conjectures and the discovery of a bijection between ASMs and
TSSCPPs.

Throughout the following we will use the notation $[x]_q$ for the usual $q$-number
\[
[x]_q = \frac{q^x-q^{-x}}{q-q^{-1}}.
\]
The notation $[x]$ will always refer to base $q$.
\section{$q$-deformed  Knizhnik--Zamolodchikov equation}

\begin{defn}
\label{def:TL} The \textbf{Temperley--Lieb algebra} of type $A_L$, denoted
by $\T^{\rm A}_L(q)$, is the unital algebra defined in terms of
generators $e_i$, $i=1,\ldots,L-1,\,$ satisfying the relations
\begin{align}
&e_i^2 = -[2] e_i, \quad e_ie_j = e_je_i \quad\forall\, i,j:\, |i-j| >1,\nonumber \\
&e_ie_{i\pm1}e_i = e_i. \label{TLdef}
\end{align}
\end{defn}
Every representation of the Temperley--Lieb algebra defines a solvable lattice model through the definition of so-called R-matrices. An R-matrix is a representation of the following Baxterised element of the Temperley--Lieb algebra:
\begin{align}
R_i(u)& =\,
\frac{[1-u] - [u]\, e_i}{[1+u]},
\end{align}
where $u\in {\mathbb C}\setminus \{-1\}$ is the spectral parameter.

\begin{defn}
\label{def:Dyck}
A \textbf{Dyck path} $\alpha=(\alpha_0,\alpha_1,\ldots,\alpha_L)$ is a sequence of integers $\alpha_i$ (heights) such that $\alpha_{i+1} = \alpha_i \pm 1$, $\alpha_i \geq 0$ and $\alpha_0=\alpha_L=0$.
\end{defn}

The Temperley--Lieb algebra has a known action on Dyck paths, which is well-documented in the literature, see for example \cite{Martin,GierP07} and references therein. The span of Dyck paths forms a module of $\T^{\rm A}_L(q)$, and we will denote its basis elements by $\ket{\alpha}$, where $\alpha$ runs over the set $\D_L$ of Dyck paths of length $L$. Let us now consider a linear combination $\ket\Psi$ of states $\ket \alpha$ with coefficients $\psi_\alpha$ taking values in the ring of formal series in $L$ variables $q^{\pm x_i}, \; i=1,2,\dots ,L$:
\[
\ket{\Psi(x_1,\ldots,x_{L})} = \sum_{\alpha\in\D_L}
\psi_\alpha(x_1,\ldots,x_{L}) \ket{\alpha}.
\]
The $q$-deformed  Knizhnik--Zamolodchikov equation
on a segment (with reflecting boundaries)
is a system of finite difference equations on the vector $\ket\Psi$.
This equation reads (see \cite{JKKMW,DF05})
\begin{align}
R_i(x_i-x_{i+1}) \ket\Psi &= \pi_i \ket\Psi,\qquad \forall\, i=1,\ldots,L-1,
\nonumber\\
\ket\Psi  &= \pi_0 \ket\Psi,
\label{$q$KZTL_TypeB2}\\
\ket\Psi &= \pi_L\ket\Psi,
\nonumber
\end{align}
where $R_i$  are the Baxterised elements of the Temperley--Lieb
algebra. The operators $R_i(x_i-x_{i+1})$ act on states $\ket{\alpha}$, whereas the
operators $\pi_i$ permute or reflect arguments of the coefficient functions:
\begin{align}
\pi_i \psi_\alpha(\ldots,x_i,x_{i+1},\ldots) &= \psi_\alpha(\ldots,x_{i+1},x_{i},\ldots),
\nonumber\\
\pi_0\psi_\alpha(x_1,\ldots)&= \psi_\alpha(-x_1,\ldots),\label{pi0}\\
\pi_L\psi_\alpha(\ldots,x_L)&= \psi_\alpha(\ldots,-\lambda-x_L).\label{piL}
\end{align}
The shift $\lambda\in {\mathbb C}$ is a parameter related to the level of the $q$KZ equation, see \cite{EFK}.

In \cite{DF07} polynomial solutions of the $q$KZ equation were studied in the limit where $x_i \rightarrow 0$. Interestingly, in this limit the coefficients $\psi_\alpha$ turn out to be polynomials with positive coefficients in $\tau^2$, where $\tau=-[2]=-q-q^{-1}$. Furthermore, based on observations on explicit solutions for small values of $L$, an intriguing connection between (generalised) sum rules and the enumeration of weighted CSTCPPs was conjectured. This was subsequently proved using multi-integral formulae \cite{DFZJ07}. Here we shall generalise this result.

\section{Explicit solutions}
\label{se:sols}

In \cite{DF07} explicit solutions of the $q$KZ equation in the
limit $x_i\rightarrow 0$ were obtained for $L\leq 8$. In order to
illustrate our results we shall here list the first few solutions.
We will write shorthand $\psi_\alpha$ for the limit
$x_i\rightarrow 0$ of $\psi_\alpha(x_1,\ldots,x_L)$.  The complete
solution is determined up to an overall normalisation, and we will
choose
\[
\psi_{\Omega} = \tau^{\lfloor {L/2}\rfloor (\lfloor L/2\rfloor -1)/2}
\]
for the coefficient corresponding to the maximal Dyck path
$\Omega \in\mathcal{D}_{L} : \Omega_i = \min\{i,L+\epsilon_L-i\}$,
$\epsilon_L=L \bmod 2$.
Together with the solution we list certain powers
$\tau^{c_{\alpha,n}}$ whose meaning will become clear below.

\begin{table}[h]
\renewcommand{\arraystretch}{1.4}
\[
\begin{array}{|c|c|lcc|}
\hline
L=4 & \alpha & \psi_\alpha & \tau^{\pm c_{\alpha,1}}  &\\ \hline
&\includegraphics[width=24pt]{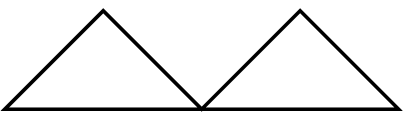} & 1+\tau^2 & 1  & \\
&\includegraphics[width=24pt]{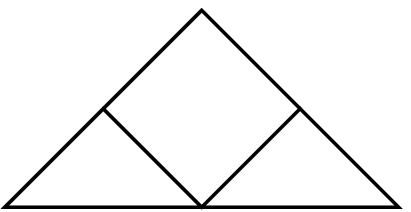} & \tau & \tau^{\pm 1} &
\\\hline\hline
L=5 & \alpha & \psi_\alpha & \tau^{\pm c_{\alpha,2}} & \tau^{\pm c_{\alpha,1}}\\ \hline
&\includegraphics[width=30pt]{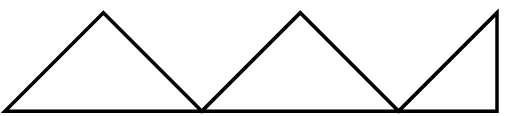} & \tau^2(2+\tau^2)& 1 & \\
&\includegraphics[width=30pt]{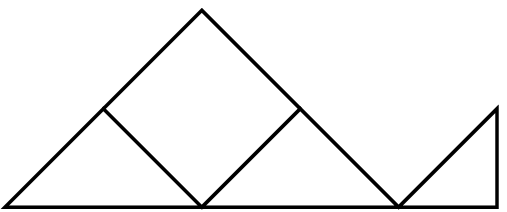} & \tau^3 & \tau^{\pm 1} & \\
&\includegraphics[width=30pt]{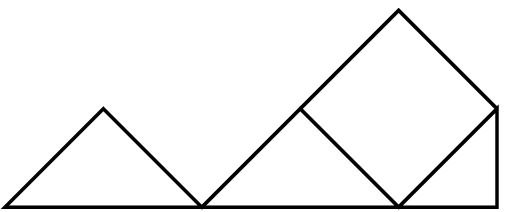} & \tau(2+\tau^2) & \tau^{\pm 1} & \\
&\includegraphics[width=30pt]{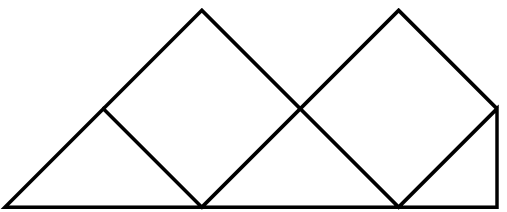} & 1+2\tau^2 & \tau^{\pm 2} & 1 \\
&\includegraphics[width=30pt]{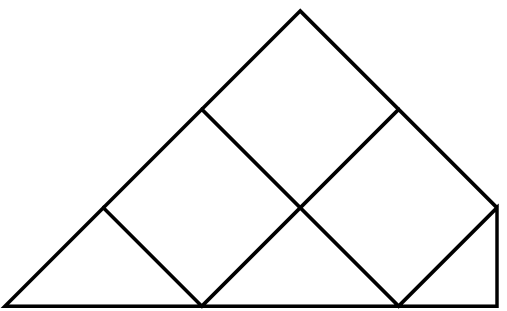} & \tau & \tau^{\pm 1} & \tau^{\pm 1}
\\\hline\hline
L=6 &\alpha & \psi_\alpha & \tau^{\pm c_{\alpha,2}} & \tau^{\pm c_{\alpha,1}} \\ \hline
&\includegraphics[width=36pt]{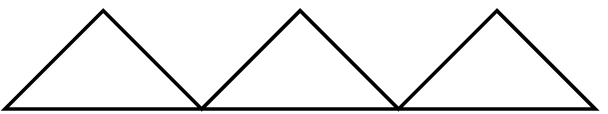} & 1+5\tau^2+4\tau^4+\tau^6 & 1&  \\
&\includegraphics[width=36pt]{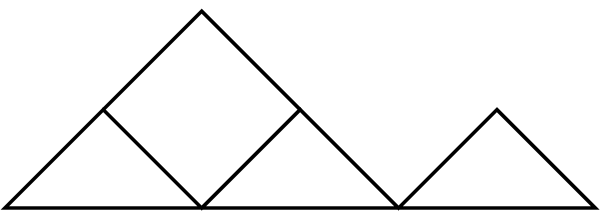} & \tau(2+2\tau^2+\tau^4) & \tau^{\pm 1} & \\
&\includegraphics[width=36pt]{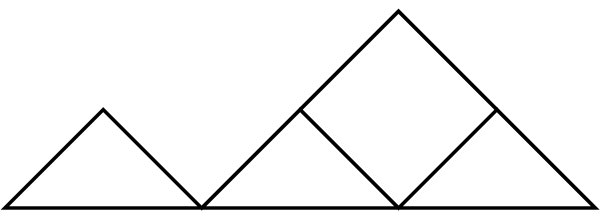} & \tau(1+3\tau^2+\tau^4) & \tau^{\pm 1}  & \\
&\includegraphics[width=36pt]{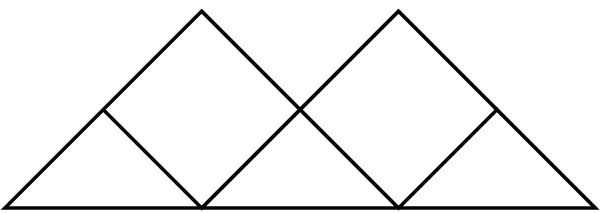} & 2\tau^2(1+\tau^2) & \tau^{\pm 2} & 1 \\
&\includegraphics[width=36pt]{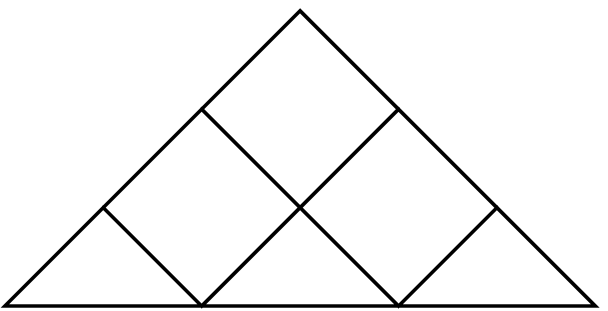} & \tau^3 & \tau^{\pm 1} & \tau^{\pm 1}\\

\hline
\end{array}
\]
\caption{Explicit solutions of the $q$KZ equation in the homogeneous limit}
\label{tab:sols}
\end{table}

\noindent
An immediate observation about these solutions was already noted in \cite{DF07}:
\begin{conj}
The components $\psi_{\alpha}(x_1,\ldots,x_L)$ of the polynomial
solution of the $q$KZ equation of type A in the limit
$x_i\rightarrow 0$, $i=1,\ldots,L$ are, up to an overall factor
which is a power of $\tau$, polynomials in $\tau^2$ with
\textbf{positive} integer coefficients. Here $\tau=-[2]$.
\end{conj}
\noindent
Polynomiality and integrality of the coefficients is proved in \cite{DFZJ07}.

\begin{figure}[h]
\centerline{\includegraphics[width=200pt]{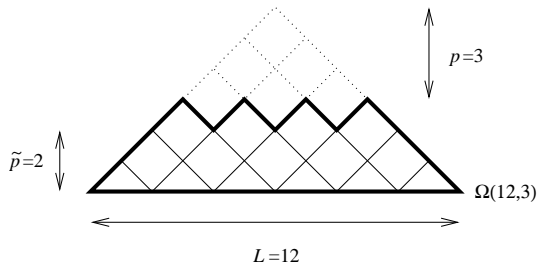}}
\caption{The  path $\Omega(12,3)\in \mathcal{D}_{12,3}$ .}
\label{fig:omegan}
\end{figure}

Based on these explicit solutions, and solutions obtained for
$L=9$ and $L=10$ in \cite{GierP07}, we discovered underlying
discrete bilinear relations. In order to uncover these relations
we have to introduce certain partial sums over the components
$\psi_\alpha$ of the solution. Let us first define the paths
$\Omega(L,p) \in \mathcal{D}_L$ whose local minima lie at height $\tilde{p}$, where
\[
\tilde{p}=\lfloor{(L-1)/2}\rfloor-p,\qquad p=0,\dots ,\lfloor{(L-1)/2}\rfloor .
\]
Figure~\ref{fig:omegan} illustrates the path $\Omega(12,3)$. We
further define the subset $\mathcal{D}_{L,p}$ of Dyck paths of
length $L$ which lie above $\Omega(L,p)$, i.e. whose local minima
lie on or above height $\tilde{p}$. Formally, this subset is
described as
\[
\mathcal{D}_{L,p} = \left\{ \alpha\in\mathcal{D}_L |\ \alpha_i\geq \Omega_i(L,p) =
\min(\Omega_i, \tilde{p}) \right\},
\]
where $\Omega_i$ are integer heights of
the maximal Dyck path $\Omega = \Omega(L,0)$.

\begin{figure}[h]
\centerline{\includegraphics[width=200pt]{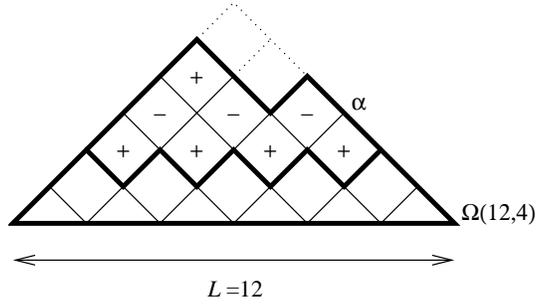}}
\caption{Definition of the number $c_{\alpha,p}$ as the signed sum of boxes between the $\alpha$ and the path $\widetilde{\Omega}(12,4)$. In this figure $L=12$ and $p=4$ and $c_{\alpha,4}=4-3+1=2$.}
\label{fig:cpmdef}
\end{figure}

To each Dyck path we associate an integer $c_{\alpha,p}$.
Let $\alpha=(\alpha_0,\alpha_1, \ldots,\alpha_L)\in
\mathcal{D}_{L,p}$ be a Dyck path of length $L$ whose minima lie
on or above height $\tilde{p}$. Then $c_{\alpha,p}$ is defined as
the signed sum of boxes between $\alpha$ and $\Omega(L,p)$, where
the boxes at height $\tilde{p}+h$ are assigned $(-1)^{h-1}$ for
$h\geq 1$. An example is given in Figure~\ref{fig:cpmdef}, and an
explicit expression for $c_{\alpha,p}$ is given by
\begin{equation}\label{cap}
c_{\alpha,p} =
\frac{(-1)^{\tilde{p}}}{2}
\sum_{i=2}^{L-1}(-1)^i(\alpha_i-\Omega_i(L,p)).
\end{equation}

In the next section we will consider certain properties of the partial weighted sums
\be
S_{\pm}(L,p) = \sum_{\alpha\in\mathcal{D}_{L,p}} \tau^{\pm c_{\alpha,p}}\, \psi_\alpha.
\label{eq:Spmdef}
\ee
The partial sums corresponding to the solutions in Table~\ref{tab:sols} are given in Table~\ref{tab:partsums}.
\begin{table}[h]
\renewcommand{\arraystretch}{1.4}
\[
\begin{array}{|c|ll|}
\hline
L=4 & S_-(4,0) = \tau & S_+(4,0) = \tau\\
& S_-(4,1) = 2+\tau^2 & S_+(4,1) = 1+2\tau^2\\
\hline\hline
L=5 & S_-(5,0) = \tau & S_+(5,0) = \tau\\
& S_-(5,1) = 2(1+\tau^2)& S_+(5,1) = 1+3\tau^2 \\
& S_-(5,2) = \tau^{-2}(1+5\tau^2 + 4\tau^4 + \tau^6) & S_+(5,2) = \tau^2(6+5\tau^2)\\
\hline\hline
L=6 & S_-(6,0) = \tau^3 & S_+(6,0) = \tau^3\\
& S_-(6,1) = \tau^2(3+2\tau^2) & S_+(6,1) = \tau^2(2+3\tau^2) \\
& S_-(6,2) = 6+13\tau^2+6\tau^4+\tau^6 & S_+(6,2) = 1+8\tau^2+12\tau^4+5\tau^6\\
\hline
\end{array}
\]
\caption{Partial sums corresponding to the solutions in Table~\ref{tab:sols}.}
\label{tab:partsums}
\end{table}

\section{Plane partitions and the discrete Hirota equation}
\label{PP+H} It was observed in \cite{Pyatov} that for $\tau=1$
($q=\e^{2\pi\i/3}$), the partial sums in Table~\ref{tab:partsums}
satisfy a discrete bilinear relation called Pascal's hexagon, or
the discrete Boussinesq equation which, in turn, is a two
dimensional reduction of the discrete Hirota equation (see
\cite{Zabrodin} and references therein). Here we generalise this
result to arbitrary $\tau$, and show that the partial sums satisfy
the discrete Hirota equation. Based on experimental data for
solutions up to $L=10$ \cite{GierP07}, we were led to introduce
the following polynomials:~\footnote{The idea to consider such
type determinants comes from the observation that some of the
polynomials given in Table~\ref{tab:partsums} coincide with the
conjectural generating functions of VSASMs  and VHSASMs
$T_n(\tau^2,\mu)|_{\mu=0, 1}$ (see \cite{Robbins}, Table 4.3). The
latter functions are given by similar determinants.}
\be
T(L,p,k) = \det_{1\leq \ell,m \leq p} \left( \sum_{r=0}^{2p} \binom{\ell+k-1}{r-\ell} \binom{m+L-2p-k}{2m-r} \tau^{2(2m-r)}\right).
\label{Tdef}
\ee
We then have the following result for the functional form of the
partial sums for arbitrary system size $L$:
\begin{prop}
\label{conj:S=T}
\begin{align}
S_+(L,p) &= \tau^{\nu_{L,p}}\, T(L,p,\lfloor L/2\rfloor -p),\\
S_-(L,p) &= \tau^{\nu_{L,p}}\, T(L,p,\lfloor L/2\rfloor -p+1),
\end{align}
where
\[
\nu_{L,p} = \frac12 \bigl(\lfloor L/2 \rfloor (\lfloor L/2 \rfloor -1) -p(p+1)\bigr).
\]
\end{prop}
\noindent
This proposition will be proved using multiple integral equations in Section~\ref{se:proof}.

\subsection{Punctured symmetric plane partitions}
\label{se:ppp}

\begin{figure}[h]
\centerline{\includegraphics[width=4cm]{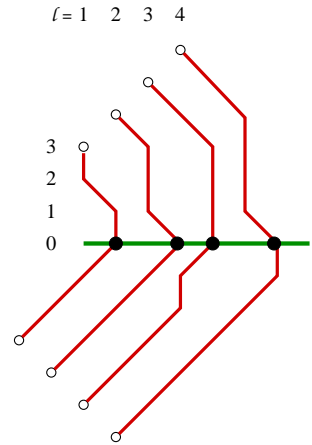}}
\caption{An example of two sets of paths with common endpoints enumerated by $T(L,p,k)$ for $L=13$, $p=4$ and $k=3$.}
\label{fig:paths}
\end{figure}

Using
the standard presentation of the determinant of $p\times p$ matrix $A$
in terms of the minors of $p\times 2p\;$ matrices $B$ and $C$: $A=B C^t$,
the polynomials
$T(L,p,k)$ can alternatively be written as
\begin{multline}
T(L,p,k) =  \sum_{1\leq r_1 <r_2 <\ldots <r_p} \det_{1\leq \ell,m \leq p} \left( \binom{\ell+k-1}{r_m-\ell} \right)
\\
\times \det_{1\leq \ell,m \leq p} \left( \binom{\ell+L-2p-k}{2\ell-r_m} \tau^{2(2\ell-r_m)} \right).
\end{multline}
By the use of Lindstr\"{o}m-Gessel-Viennot formula \cite{LGV}
the last expression can be interpreted as the generating
function of two sets of paths with the same end-point, as in
Figure~\ref{fig:paths}. One set of paths starts at positions
$(\ell,\ell+k-1)$, ($\ell=1,\ldots, p$) and has diagonal NW-SE and vertical steps; the other set starts at
$(\ell-L+2p+k,-\ell-L+2p+k)$ with diagonal NE-SW and vertical steps. Each vertical step below the horizontal line (green online) is assigned a weight $\tau^2$. All other steps are assigned weight 1. Up to an overall factor $\tau^{(k-1)p}$ this is also the generating function of pairs of paths where all vertical steps are assigned a weight $\tau$.

\begin{figure}[h]
\centerline{\includegraphics[width=10cm]{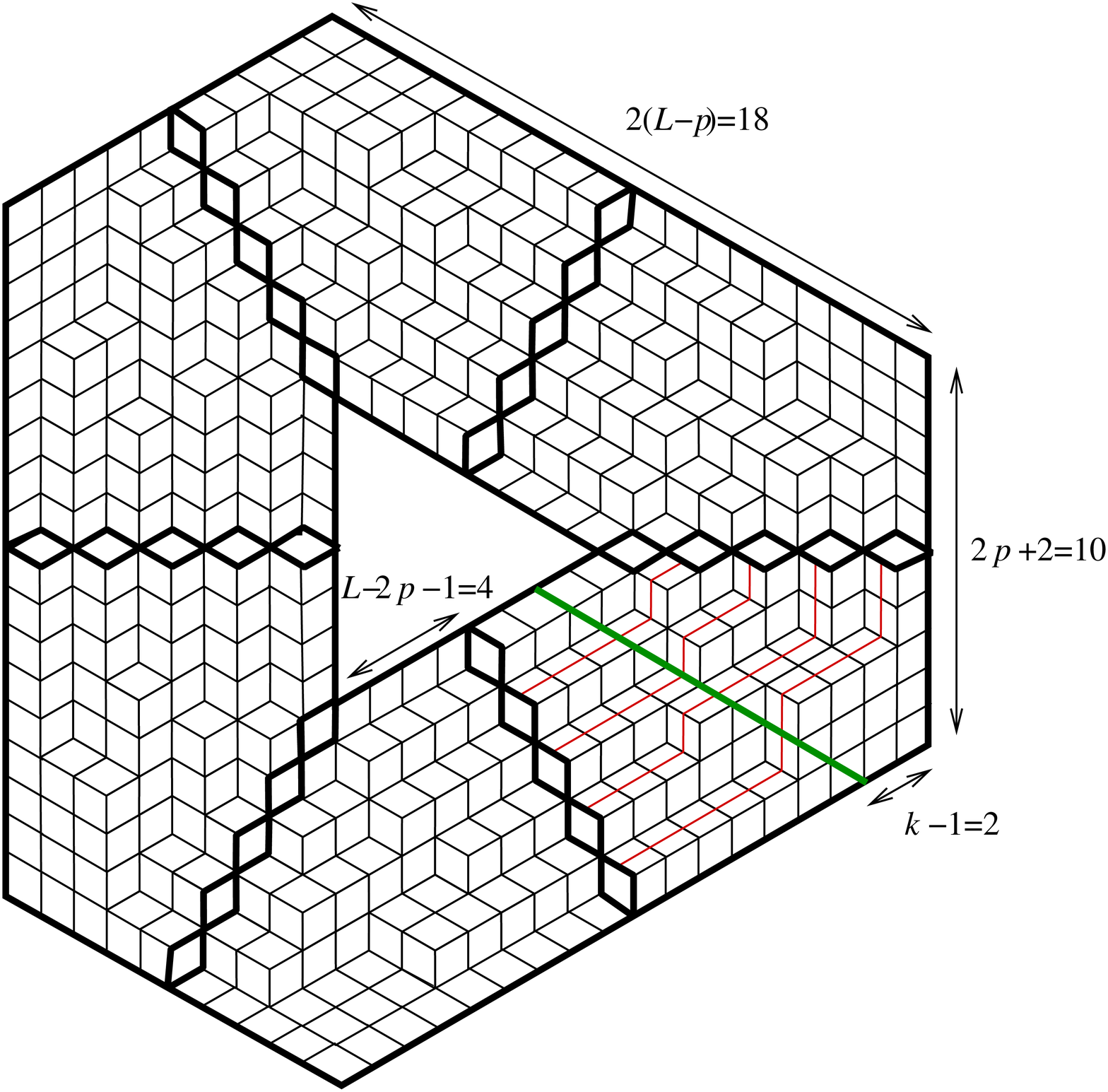}}
\caption{A punctured cyclically symmetric transpose complement
plane partition (PCSTCPP) and its fundamental region. The linear
size of the puncture is given by the difference of the lengths of
the sides of the hexagon. The position of the bisecting line
(green online) determines a particular weighting of the PCSTCPP:
vertical steps below this line are assigned a weight $\tau^2$.}
\label{fig:PCSTCPP}
\end{figure}

Consider stacking a large cube with smaller cubes, starting in a
corner of the large cube, in such a way that the smaller cubes can
only be stacked on top or next to each other. Such stackings are
called plane partitions, which are well known to be equivalent to
rhombus tilings of a hexagon. The paths in Figure~\ref{fig:paths}
arise naturally when one considers cyclically symmetric transpose
complement plane partitions (CSTCPPs) and generalisations thereof,
as was done by Ciucu and Krattenthaler \cite{CK}. CSTCPPs
correspond to rhombus tilings of a hexagon which are invariant
under rotations over $2\pi/3$ as well as under reflections across
a symmetry axis not passing through the hexagon corners. In
\cite{CK} the more general problem was considered of cyclically
symmetric transpose complement tilings of a hexagon with a
triangular hole, as in Figure~\ref{fig:PCSTCPP}. We will call such
tilings punctured CSTCPPs, or PCSTCPPs, the size of the triangular
puncture being determined by the difference of the lengths of the
sides of the hexagon.

If we weight PCSTCPPs by assigning a weight $\tau^2$ to each
vertical step below the bisecting line (green online) in the
South-East region of Figure~\ref{fig:PCSTCPP}, the weighted
enumeration of tilings of the fundamental domain of PCSTCPPs is
equivalent to the weighted enumeration of paths in
Figure~\ref{fig:paths}. We have sketched the bijection in
Figure~\ref{fig:PCSTCPP}, full details may be found in \cite{CK}
where the enumeration ($\tau=1$) is considered, which is shown to
factorise completely, see also \eqref{Sfactor}. For general $\tau$
but $p$ restricted to $p=\lfloor (L-1)/2\rfloor$ PCSTCPPs were
considered in the context of the $q$KZ equation by Di Francesco
\cite{DF07}.

\begin{figure}[h]
\centerline{\includegraphics[width=6.5cm]{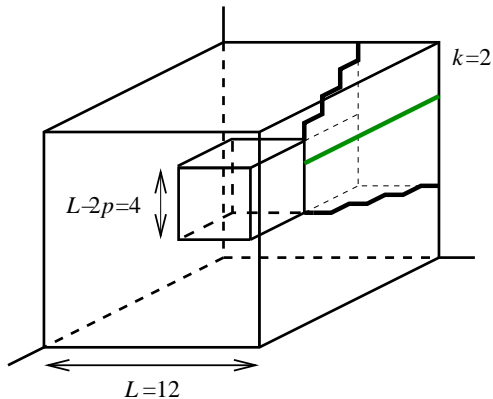}}
\caption{Fundamental region for a subclass of punctured cyclically
symmetric self complement plane partitions (PCSSCPPs). The
position of the bisecting line (green online) determines a
particular weighting of a \mbox{PCSSCPP.} } \label{fig:box}
\end{figure}

Interestingly, there is another way to interpret the paths in
Figure~\ref{fig:paths} in terms of punctured cyclically symmetric
plane partitions. When we remove a central cubic region from the
large cube, and consider special sixfold rotational symmetric
rhombus tilings of the corresponding punctured hexagon by defining
a fundamental region as in the upper right corner of
Figure~\ref{fig:box}, we obtain a subset of punctured cyclically
symmetric self complement plane partitions (PCSSCPPs).

\begin{figure}[h]
\centerline{\includegraphics[width=10cm]{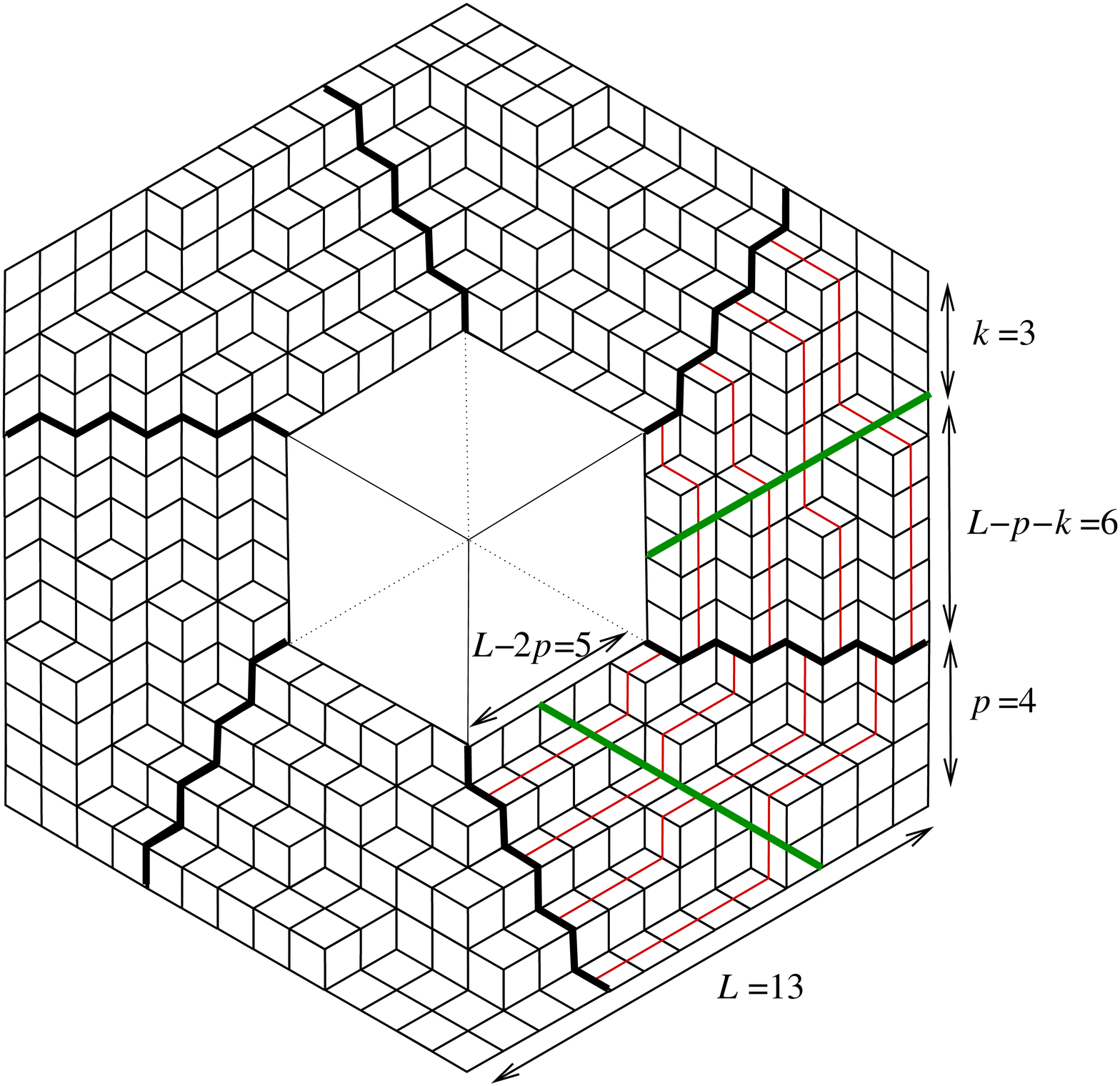}}
\caption{A PCSSCPP of size $L=13$ with a puncture of size
$L-2p=5$. PCSSCPPs are enumerated by non-intersecting lattice
paths, as indicated by the red lines. The paths in this picture
define a natural subset of all PCSSCPPs, and correspond to those
in Figure~\ref{fig:paths}. In the East region, diagonal steps
below the bisecting line (green online) are assigned a weight
$\tau^2$. All other steps are assigned weight 1. The value of $k$,
$k=3$ in this figure, determines the position of the bisecting
green line.} \label{fig:tiling}
\end{figure}

PCSSCPPs are enumerated by giving a special weight to the
non-intersecting lattice paths in Figure~\ref{fig:paths}, see \cite{CK}.
Here we will not discuss this, but will only consider
those PCSSCPPs as defined in Figure~\ref{fig:tiling}. The paths in
this picture correspond to those in Figure~\ref{fig:paths}. In the
East region of Figure~\ref{fig:tiling}, diagonal steps below the
bisecting line (green online) are assigned a weight $\tau^2$.

\subsection{Discrete Hirota equation}
We end this section with the following important observation for which we do not know a good interpretation. The polynomials $T(L,p,k)$ defined in \eqref{Tdef} satisfy a discrete bilinear relation:
\begin{prop}
\label{prop-Trecur}
The polynomials \eqref{Tdef} satisfy the recurrence
\begin{multline}
T(L,p,k)T(L-2,p-2,k+2)=T(L-1,p-2,k+2) T(L-1,p,k)\\+ \tau^2T(L-2,p-1,k)T(L,p-1,k+2).
\label{Trecur}
\end{multline}
\end{prop}
\noindent
This proposition will be proved in Section~\ref{HirotaProof}.

The recurrence \eqref{Trecur} is equivalent to the equation with which we
began this paper: \eqref{Trecur} becomes the discrete Hirota equation \eqref{Hirota}
or octahedron recurrence, when we change variables to $n=L-p-k$, $i=2p+k-L$, $j=p+k$,
and $T(L,p,k)=f(n,i,j)$:
\begin{multline*}
f(n,i,j)f(n-2,i,j)=f(n-1,i-1,j)f(n-1,i+1,j)\\+\tau^2f(n-1,i,j-1)f(n-1,i,j+1).
\end{multline*}

\section{Fully packed loop diagrams and the Razumov--Stroganov conjecture}
\label{se:fpl}

In this section we will discuss another combinatorial
interpretation of the solutions of the $q$KZ equation. At $\tau=1$
the $q$KZ equation is equivalent to an eigenvalue equation for the
transfer matrix of the O(1) loop model. In this context a relation
was conjectured in \cite{PRGN} (as a variant of a similar
conjecture by Razumov and Stroganov \cite{RS1,BGN01,RS2}) between the
solutions in Section~\ref{se:sols} for $L$ even and refined enumerations of
vertically symmetric alternating-sign matrices (VSASMs). For $L$ odd there is a connection to related objects, see below.
The $q$KZ equation is a generalisation of the O(1) eigenvalue equation to $\tau\neq 1$, as
found by Pasquier \cite{Pasquier}, and Di Francesco and Zinn-Justin \cite{DFZJ05}.

The Razumov--Stroganov conjecture can be roughly described as
follows. First a simple bijection of VSASMs to fully packed loop diagrams is made, see
e.g. \cite{Propp}. To each such fully packed loop diagram is associated a Dyck path; this will be described below. For a given Dyck path $\alpha$ there correspond many FPL diagrams, and their
number is precisely $\psi_\alpha$ at $\tau=1$. For general $\tau$
there is as yet no interpretation of $\psi_{\alpha}$ as a weighted
enumeration of FPL diagrams or VSASMs. However, as we will show in Proposition~\ref{PP2ASM} at the end of this section, there is such an interpretation for the total $\tau$-normalisations $S_{-}(2n,n-1)$ of the $q$KZ solution. We will see that $S_{+}(2n,n-1)$ and $S_{-}(2n-1,n-1)$ are also related to ASM generating functions.

As shown in the previous section, the normalisation $S_-(2n,n-1)$ is also equal to a generating function of $\tau^2$-weighted PCSTCPPs with a small puncture of size 2, establishing a connection between such plane partitions and weighted VSASMs. A similar connection was considered in \cite{DF07}. Other equalities between generating functions of weighted cyclically symmetric plane partitions and generating functions related to weighted enumerations of symmetric ASMs were conjectured by Robbins and Kuperberg \cite{Robbins,Kuperberg}. We will prove some of these in Section~\ref{ProofofPP2ASM}.

Finding the $\tau$-statistic on FPL diagrams that gives rise to an interpretation
of $\psi_\alpha$ for $\tau\neq 1$  should provide clues for explicit bijections between (P)CSTCPPs and symmetry classes of ASMs, as well as for an explicit bijection between TSSCPPs and unrestricted ASMs. In fact, in this section we will make another small step towards such a bijection by elucidating the role played by the parameter $p$ which determines the size of the puncture in PCSTCPPs.

\subsection{Fully packed loop diagrams}
\begin{figure}[h]
\centerline{\includegraphics[width=9cm]{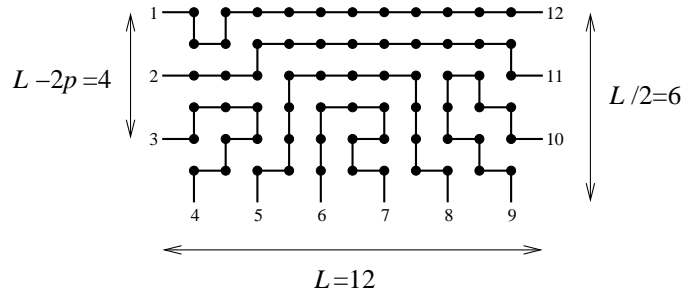}}
\caption{An FPL diagram of size $L=12$ with $L/2 -p =2$ loop lines connecting
loop terminals $1,2$ with $L,L-1$ respectively.}
\label{fig:fpl_even}
\end{figure}

\begin{figure}[h]
\centerline{\includegraphics[width=10cm]{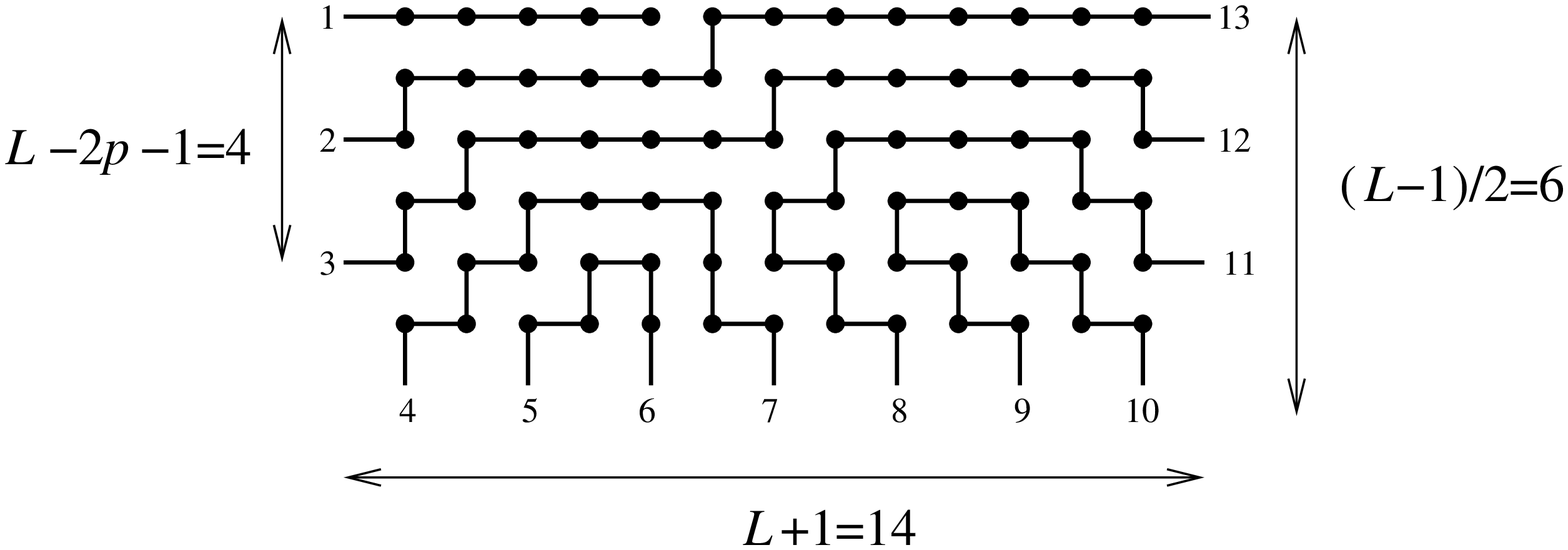}}
\caption{An FPL diagram of size $L=13$ with $\lfloor L/2\rfloor -p =2$ loop
lines connecting external loop terminals $2,3$ with $L,L-1$ respectively.}
\label{fig:fpl_odd}
\end{figure}

For $L$ even, consider the set of vertices of a piece of square
lattice of size $L/2 \times L$. In the case of $L$ odd, take a
piece of size $(L-1)/2 \times (L+1)$. On this set of points, draw
bonds such that each internal vertex has exactly two drawn bonds.
On the boundary we impose the condition that every other ingoing
bond on the left hand side, bottom and right hand side is drawn,
starting at the top left vertex. \footnote{These boundary
conditions and the shape of the lattice are of particular
relevance in this paper, but more general FPL diagrams can be
considered, see e.g. \cite{DeGier}.} In this way, the drawn bonds
form closed polygons or connect outgoing bonds to each other, see
for example Figures~\ref{fig:fpl_even} and \ref{fig:fpl_odd}. From
Kuperberg's work \cite{Kuperberg} it follows that FPL diagrams of
even size are equinumerous with VSASMs. Those of odd size are
conjectured\footnote{To our knowledge no proof of this assertion
exists.} to be equinumerous with CSTCPPs \cite{PRGN}.

\begin{figure}[h]
\centerline{\includegraphics[width=4cm]{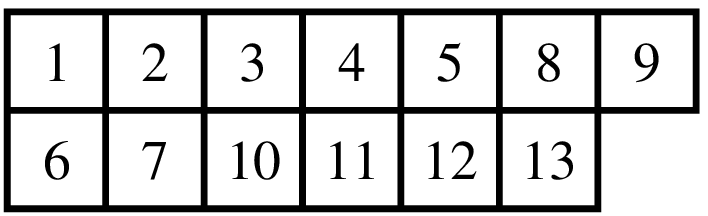}}
\caption{Standard Young tableau corresponding to the FPL diagram in
Figure~\ref{fig:fpl_odd}.}
\label{fig:young}
\end{figure}

We will label the outgoing bonds by successive integers, as in
Figures~\ref{fig:fpl_even} and \ref{fig:fpl_odd}. As each outgoing
bond is connected to another outgoing bond, FPL diagrams can be
naturally labeled by link patterns, or equivalently, two-row Young
tableaux or Dyck paths. For example, the diagram in
Figure~\ref{fig:fpl_odd} has link pattern $(\, ((\, (())\, (())\,
))$ which is short hand for saying that $1$ is connected to the
top of the diagram, $2$ is connected to $13$, $3$ to $12$ etc. In
general, to each link pattern correspond many FPL diagrams. This
information can also be coded in two-row standard Young tableaux.
The entries of the first row of the Young tableau correspond to
the positions of opening parentheses '(' in a link pattern, and
the entries of the second row to the positions of the closing
parentheses ')'. The FPL diagram of Figure~\ref{fig:fpl_odd}
carries as a label the standard Young tableau given in
Figure~\ref{fig:young}.

\begin{figure}[h]
\centerline{\includegraphics[width=5cm]{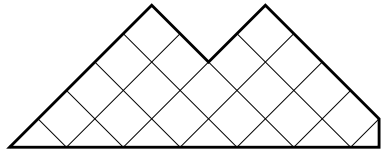}}
\caption{Dyck path corresponding to the FPL diagram in Figure~\ref{fig:fpl_odd}
and the standard Young tableau in Figure~\ref{fig:fpl_odd}.}
\label{fig:dyck}
\end{figure}

Yet another way of coding the same information is by using Dyck
paths. This will be useful when making a connection to the results
of Section~\ref{se:sols}. Each entry in the first row of the
standard Young tableau represents an up step, while those in the
second row represent down steps. The Dyck path corresponding to
Figure~\ref{fig:young} is given in Figure~\ref{fig:dyck}.

\subsection{Subsets of FPL diagrams and punctured plane partitions}

As discussed above, FPL diagrams of even size are equinumerous
with VSASMs, and it follows from the results of Section~\ref{PP+H}
that they are also equinumerous to PCSTCPPs with a small puncture
of size 2. Here we shall formulate a refined correspondence
between subsets of FPL diagrams and PCSTCPPs.

Recall the subset $\D_{L,p}$ of Dyck paths whose local minima lie
on or above height $\tilde{p}$ where $\tilde{p}=\lfloor
(L-1)/2\rfloor -p$, see Section~\ref{se:sols}. Each path in
$\D_{L,p}$ is a label for FPL diagrams whose loop terminals
$1,\ldots \tilde{p}$ (in the case of $L$ even) are connected to
terminals $L-\tilde{p}+1,\ldots,L$ respectively, see
Figure~\ref{fig:fpl_even}. In the case of $L$ odd, the Dyck paths
in $\D_{L,p}$ label FPL diagrams with loop connections between
terminals $2,\ldots,\tilde{p}+1$ and $L-\tilde{p}+1,\ldots,L$
respectively, see Figure~\ref{fig:fpl_odd}.

\begin{defn}
A \textbf{$\mathbf{\textit{p}}$-restricted FPL diagram} is an FPL
diagram whose corresponding Dyck path belongs to $\D_{L,p}$.
\end{defn}
\noindent In other words, in a $p$-restricted FPL diagram of even
size, the first $L/2-p$ loop terminals are connected to the last
$L/2-p$ loop terminals. This classification of FPL diagrams allows
us to formulate a refined correspondence between FPL diagrams and
symmetric plane partitions:

\begin{conj}
\label{conj:FPL2PPP} The total number of $p$-restricted FPL
diagrams of size $L=2n$ is equal to the total number $\left.
S_{\pm}(2n,p)\right|_{\tau=1}$ of PCSTCPPs having sides of length $2(L-p)$ and $2(p+1)$ with a
triangular puncture of size $2(L-2p-1)$, see Figure \ref{fig:PCSTCPP}.
\end{conj}
\noindent As a side remark we note that $\left.
S_{\pm}(L,p)\right|_{\tau=1}$ at $\tau=1$ factorises (this was
proved in \cite{MillsRR,CK} in the context of PCSTCPPs) and takes
the form
\be \left. S_{\pm}(L,p)\right|_{\tau=1} = 2^{(p+1)(L-p)}
\prod_{j=1}^p \frac{\Gamma(L-j+1) \Gamma((2L+2j+3)/6)
\Gamma((L-2j+3)/3)}
{\Gamma(L-2j+1)\Gamma(j+1/2)\Gamma((2L-j+3)/3)}.
\label{Sfactor}
\ee
This expression is equivalent to earlier conjectured expressions
for the total number of $p$-restricted FPL diagrams in the context
of the O(1) loop model and the Razumov-Stroganov conjecture
\cite{MNGB,Pyatov}.

One may ask if a generalisation of Conjecture~\ref{conj:FPL2PPP}
holds for $\tau\neq 1$. As we have seen in Section~\ref{se:ppp},
the parameter $\tau^2$ is a natural weight for (punctured)
symmetric plane partitions. One can therefore ask whether it also
describes a natural statistic on ASMs or FPL diagrams. In the case
of $p=\lfloor (L-1)/2\rfloor$ ($\tilde{p}=0$) this question can be
answered in the affirmative. It was already observed by Robbins and
Kuperberg \cite{Robbins,Kuperberg} that some $\tau^2$-generating
functions for ordinary CSTCPPs, i.e. without puncture, are the
same as certain polynomials arising in $\tau^2$-enumerations of
symmetric ASMs where each $-1$ is assigned a weight $\tau^2$. In
the case of FPL diagrams this amounts to giving a weight $\tau$ to
every two consecutive vertical or horizontal steps (such
consecutive steps correspond to either a $+1$ or $-1$ in ASM
language). Our expressions for the partial sums in
Table~\ref{tab:partsums} can be directly compared with Kuperberg's
Table 4. Using the notations of \cite{Kuperberg}, see also
Section~\ref{ProofofPP2ASM}, we collect this observation with two
others in the following intriguing proposition:

\begin{prop}
\label{PP2ASM}
\begin{align*}
T(2n,n-1,2) &= S_-(2n,n-1) = A_{\rm V}(2n+1;\tau^2),\\
T(2n,n-1,1) &= S_+(2n,n-1) = A_{\rm VHP}^{(2)}(4n+2;\tau^2),\\
T(2n-1,n-1,1) &= \tau^{n-1} S_-(2n-1,n-1) = \tilde{A}_{\rm UU}^{(2)}(4n;\tau^2).
\end{align*}
\end{prop}
This proposition will be proved in Section~\ref{ProofofPP2ASM}.

It is an open problem to find a generalisation of
Conjecture~\ref{conj:FPL2PPP} to arbitrary weight $\tau$, or a
generalisation of Proposition~\ref{PP2ASM} to arbitrary $p$, i.e.
to find a correspondence between $\tau^2$-weighted $p$-restricted
FPL diagrams and $\tau^2$-enumerations of PCSTCPPs.

\section{Proof of the determinant formula for partial sums}
\label{se:proof}
In this section we prove Proposition~\ref{conj:S=T} using the formalism developed
in \cite{DFZJ07,DFZJ07b} and which consists of writing integral formulae for
solutions of the $q$KZ equation. For the sake of simplicity, we shall work out
separately the two possible parities of the size $L$.
\subsection{Even size}
Assume $L=2n$. The following set of integrals
was introduced in \cite{DFZJ07}:
\begin{multline}\label{integinhom}
\psi_{a_1,\ldots,a_n}=
\prod_{1\leq i<j\leq L}
[1+x_i-x_j][1-x_i-x_j]
\oint\cdots\oint \prod_{\ell=1}^n
{\log q\over q-q^{-1}} {\d y_\ell \over \pi \i}\\
{\prod_{1\leq \ell<m\leq n}
[y_\ell-y_m][1+y_\ell-y_m][y_\ell+y_m]
\prod_{1\leq \ell\leq m\leq n} [1+y_\ell+y_m]
\over \prod_{\ell=1}^n\prod_{i=1}^{L} [y_\ell+x_i]
\prod_{i=1}^{a_\ell}[y_\ell-x_i]\prod_{i=a_\ell+1}^{L}
[1+y_\ell-x_i]
}
\end{multline}
where the contour integrals surround the poles at $x_i-1$.
Here $a_1,\ldots,a_n$ form a non-decreasing sequence of integers between $1$ and $L-1$.

The relation to the solution of the $q$KZ system
(\ref{$q$KZTL_TypeB2}) is as follows: up to normalization by a symmetric
factor of the parameters $x_i$, the $\psi_{a_1,\ldots,a_n}$ are linear
combinations of the components of the solution:
\begin{equation}
\psi_{a_1,\ldots,a_n}=\sum_\alpha C_{a_1,\ldots,a_n; \alpha} \psi_\alpha
\end{equation}
with coefficients $C_{a_1,\ldots,a_n; \alpha}$ that are described explicitly in appendix
A of \cite{DFZJ07b}, and which we shall define here by recurrence.
First, $C_{\emptyset;\emptyset}=1$.
Next, for a pair $(a_1,\ldots,a_n;\alpha)$ of length $L$,
consider any local maximum $i$ of the path $\alpha$, and the new path $\alpha'$
obtained by removing the two steps before and after $i$. Call $k$ the number
of $\ell$ such that $a_\ell=i$.
If $k=0$, $C_{a_1,\ldots,a_n;\alpha}=0$.
If $k>0$, consider the new sequence $a'_1,\ldots,a'_{n-1}$ obtained from
$a_1,\ldots,a_n$ by removing one ``$i$'' and replacing each remaining $a_\ell$ with:
itself if $a_\ell<i$; $i-1$ if $a_\ell=i$; $a_\ell-2$ if $a_\ell>i$.
Then $C_{a_1,\ldots,a_n;\alpha}=[k] C_{a'_1,\ldots,a'_{n-1};\alpha'}$
(and this definition is independent of the choice of local maximum).

Since we are interested in the values of the $\psi_\alpha$ at $x_i=0$,
let us set $x_i=0$ in (\ref{integinhom}) and perform the change of variables
$u_\ell=[1+y_{n+1-\ell}]/[y_{n+1-\ell}]$. At this stage it is convenient
to reindex the integers as $b_\ell=L-a_{n+1-\ell}$, $\ell=1,\ldots,n$, and to define
$\bar\psi_{b_1,\ldots,b_n}:=\psi_{a_1,\ldots,a_n}$, so that
\begin{multline}\label{integhom}
\psi_{a_1,\ldots,a_n}=\bar\psi_{b_1,\ldots,b_n}
= \oint\cdots\oint \prod_{\ell=1}^n {\d u_\ell\over 2\pi\i\, u_\ell^{b_\ell}} \left[
\prod_{1\leq \ell \leq m \leq n} (1-u_\ell u_m) \right.\\
\times \left. \prod_{1\leq \ell < m \leq n} (u_m-u_\ell)
(1+\tau u_m+u_\ell u_m) (\tau+u_\ell+u_m) \right]
\end{multline}
Here the normalization is chosen in such a way that
$\psi_{1,2,\ldots,n}=\psi_\Omega=\tau^{n(n-1)/2}$ (note that in this case the integrals
are trivial and can be performed by simply setting $u_\ell=0$ in the numerator of the integrand).

We now consider specific $\bar\psi_{b_1,\ldots,b_n}$ which will reproduce
our sums $S_\pm(L,p)$. Fix a non-negative integer $p$, and let
$\tilde{p}=n-1-p$.
Consider sequences $(b_1,\ldots,b_n)$ of the form
\[
b_\ell=
\begin{cases}
\ell & 1\le\ell\le \tilde p+1\\
2\ell-\tilde p-1-\epsilon_{\ell-\tilde p-1}& \tilde p+2\le\ell\le n
\end{cases}
\]
where $\epsilon_1,\ldots,\epsilon_p\in \{ 0,1 \}$. We have the following
\begin{lemma}
\[
\bar\psi_{b_1,\ldots,b_n}=
\bar\psi_{1,\ldots,\tilde p+1,
\tilde p+3-\epsilon_1,\ldots,2n-\tilde p-1-\epsilon_p}
=\sum_{
\scriptstyle\alpha\in\mathcal{D}_{L,p}
\atop
\scriptstyle\forall \ell,\alpha_{L-\tilde p-2\ell-1}<\alpha_{L-\tilde p-2\ell}\
{\rm iff}\ \epsilon_\ell=1}
\psi_\alpha
\]
\end{lemma}
\begin{proof}
We shall proceed by induction. Fix a sequence of integers $b_\ell$
as in the lemma, that is in terms of the mirror-symmetric sequence $a_\ell$,
\[
a_\ell=
\begin{cases}
2\ell+\tilde p-1+\epsilon_{p+1-\ell}& 1\le\ell\le p\\
\ell+n-1& p+1\le\ell\le n
\end{cases}
\]
Note that $a_\ell\ge\tilde p+1$.
Let $\alpha$ be a Dyck path. Consider a local maximum $i$ of $\alpha$.
One can always assume $i\le n$.
There are two cases:

1. $i\le \tilde p$. In this case, $\alpha\not\in\mathcal{D}_{L,p}$.
We find immediately that there are zero $a_\ell=i$, so that the coefficient is zero.

2. $\tilde p< i \le n< L-\tilde p$. Call $\ell=\lfloor (i-\tilde p+1)/2 \rfloor$.
There are four cases depending on the parity of $i$ and the value of
$\epsilon_{p+1-\ell}$.
If $\epsilon_{p+1-\ell}=0$ and $i=2\ell+\tilde p\ne a_\ell$,
there are no $a_\ell$ equal to $i$ so that the coefficient is zero.
Since $i$ is a local maximum, we have
$\alpha_{L-\tilde p-2(p+1-\ell)-1=i-1}<\alpha_{L-\tilde p-2(p+1-\ell)=i}$ satisfying the inequality
in the summation of lemma 1 (despite $\epsilon_{p+1-\ell}=0$).
Similarly, if $\epsilon_{p+1-\ell}=1$ and $i=2\ell+\tilde p-1\ne a_\ell$,
there are no $a_\ell$ equal to $i$ so that the coefficient is zero,
and $\alpha_{L-\tilde p-2(p+1-\ell)-1=i}>\alpha_{L-\tilde p-2(p+1-\ell)=i+1}$ violating the inequality
in the summation of lemma 1 (despite $\epsilon_{p+1-\ell}=1$).
In the other two cases, we have $a_\ell=i$ and the equivalence in the summation
of lemma 1 is valid.
We can then apply the definition by recurrence of
the coefficient $C_{a_1,\ldots,a_n;\alpha}$; the new sequence $a'_1,\ldots,a'_{n-1}$
is exactly the same type as $a_1,\ldots,a_n$, that is
defined by the same $\epsilon_i$ with $\epsilon_{p+1-\ell}$ skipped.
On the other hand it is clear that the other conditions on $\alpha$ in the summation
of lemma 1 are equivalent to the conditions on $\alpha'$ ($\alpha$
with the local maximum removed)
with the new sequence $a'_1,\ldots,a'_{n-1}$.
One then uses the induction hypothesis to conclude.
\end{proof}

\begin{example}
Consider the two sequences with $L=6$, $p=1$ ($\tilde p=1$).
We find
\begin{align*}
\bar\psi_{1,2,4}&=\psi_{\includegraphics[width=36pt]{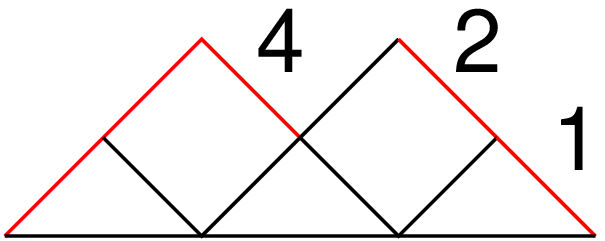}}=2\tau^2(1+\tau^2)\\
\bar\psi_{1,2,3}&=\psi_{\includegraphics[width=36pt]{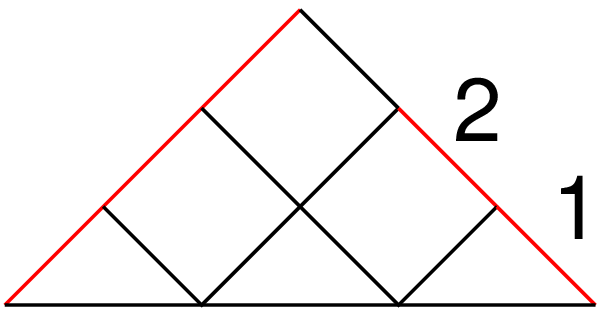}}=\tau^3\\
\end{align*}
In size $L=8$, with $p=2$ ($\tilde p=1$),
\begin{align*}
\bar\psi_{1,2,4,6}&=\psi_{\includegraphics[width=48pt]{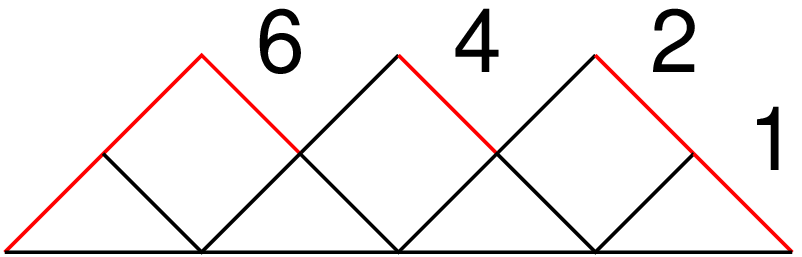}}=6\tau^3+21\tau^5+18\tau^7+5\tau^9\\
\bar\psi_{1,2,4,5}&=\psi_{\includegraphics[width=48pt]{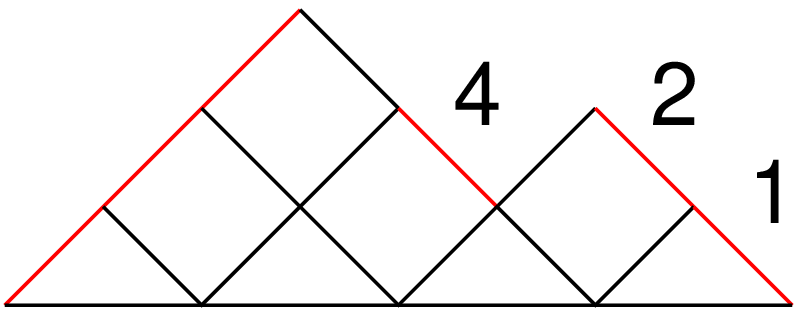}}+\psi_{\includegraphics[width=48pt]{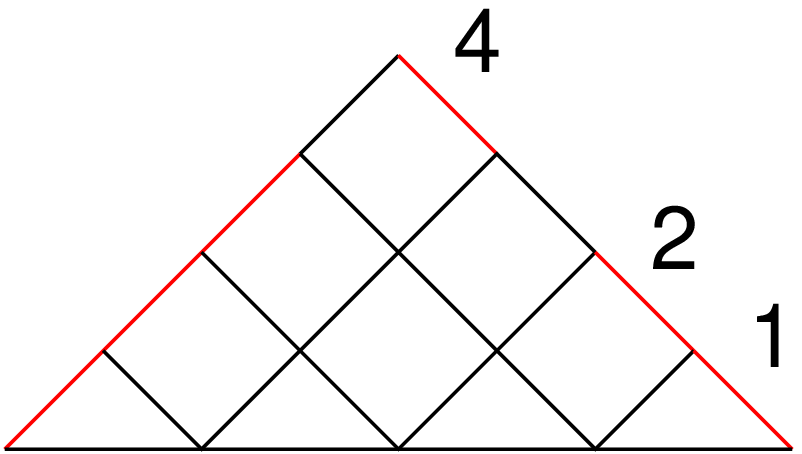}}
=5\tau^4+7\tau^6+3\tau^8\\
\bar\psi_{1,2,3,6}&=\psi_{\includegraphics[width=48pt]{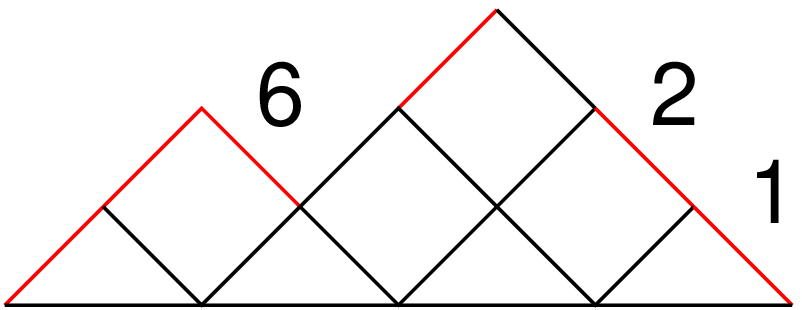}}=3\tau^4+8\tau^6+3\tau^8\\
\bar\psi_{1,2,3,5}&=\psi_{\includegraphics[width=48pt]{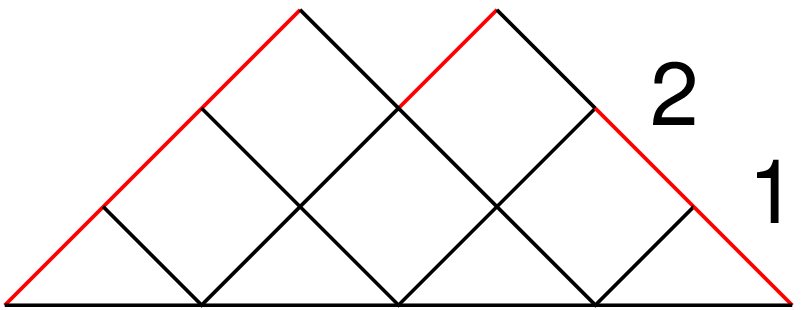}}=3\tau^5+3\tau^7\\
\end{align*}
The edges in red are those whose labels (counted from right
to left) appear in the sequence of integers iff the edge is a down step
(in fact the first and last $\tilde p+1$ edges are fixed by the fact
that the paths are in $\mathcal{D}_{L,p}$).
\end{example}

Note that taken together, the various sequences for a given $L$ and $p$ reproduce
the full set of paths of $\mathcal{D}_{L,p}$. Furthermore,
by direct computation using formula (\ref{cap}) (grouping together
pairs $\alpha_{L-\tilde p-2i-1}$ and $\alpha_{L-\tilde p-2i}$ in the sum
and using $\alpha_{L-\tilde p-2i}-\alpha_{L-\tilde p-2i-1}=2\epsilon_i-1$), it is easy to
check that all $\psi_\alpha$ that contribute to a given
$\bar\psi_{1,\ldots,\tilde p+1,\tilde p+3-\epsilon_1,\ldots,2n-\tilde p-1-\epsilon_p}$
have the same integer $c_{\alpha,p}=\sum_{i=1}^p \epsilon_i$.

We thus define
\begin{equation}
S(L,p|t)=\sum_{\epsilon_1,\ldots,\epsilon_p\in \{0,1\}}
t^{\sum_{i=1}^p \epsilon_i}
\bar\psi_{1,\ldots,\tilde p+1,
\tilde p+3-\epsilon_1,\ldots,2n-\tilde p-1-\epsilon_p}
\end{equation}
and claim that $S_\pm(L,p)=S(L,p|\tau^{\pm 1})$.

Using (\ref{integhom}), we now
obtain the following integral representation for $S(L,p|t)$:
\begin{multline}
S(L,p|t)
= \oint\cdots\oint \left( \prod_{\ell=1}^{\tilde p+1} {\d u_\ell\over 2\pi \i\, u_\ell^\ell}\right)
\left( \prod_{\ell=\tilde{p}+2}^n {\d u_\ell (1+t u_\ell)\over 2\pi \i\, u_\ell^{2\ell-\tilde{p}-1}}\right)
\left[ \prod_{1\leq \ell \leq m \leq n} (1-u_\ell u_m) \right.\\
\times \left. \prod_{1\leq \ell < m \leq n} (u_m-u_\ell)
(1+\tau u_m+u_\ell u_m) (\tau+u_\ell+u_m) \right]
\end{multline}
The first integrals over $u_1,\ldots,u_{\tilde{p}+1}$ can be performed successively by simply setting the corresponding variables to zero.
The result, after shifting the indices of the variables, is:
\begin{multline}
S(L,p|t)
= \oint\cdots\oint
\prod_{\ell=1}^p {\d u_\ell (1+t u_\ell)(1+\tau u_\ell)^{\tilde{p}+1} (\tau+u_\ell)^{\tilde{p}+1}\over 2\pi \i\, u_\ell^{2\ell}}
\left[ \prod_{1\leq \ell \leq m \leq p} (1-u_\ell u_m) \right.\\
\times \left. \prod_{1\leq \ell < m \leq p} (u_m-u_\ell)
(1+\tau u_m+u_\ell u_m) (\tau+u_\ell+u_m) \right]
\end{multline}
Next we use the following lemma:
\begin{lemma}
If $AS$ designates antisymmetrisation:
$AS(f(u_1,\ldots,u_p))=\sum_{\sigma\in \mathcal{S}_{p}} (-1)^\sigma f(u_{\sigma(1)},\ldots,u_{\sigma(p)})$, and $(\cdots)_{\le 0}$ means keeping only non-positive powers of a Laurent
polynomial in the variables
$u_\ell$, then the following equality holds:
\begin{multline}
\left\{\prod_{1\le\ell\le m\le p}(1-u_\ell u_m)\
{\rm AS}\left(\prod_{\ell=1}^p u_\ell^{-2\ell+1} \prod_{1\le\ell<m\le p}
(1+u_\ell u_m+\tau u_m)\right)\right\}_{\le 0}\\
={\rm AS}\left(\prod_{\ell=1}^p
\left(u_\ell^{-\ell}(\tau +u_\ell^{-1})^{\ell-1}\right)\right)=
\prod_{\ell=1}^p u_\ell^{-1}\prod_{1\le\ell<m\le p}(u_m^{-1}-u_\ell^{-1})
(\tau+u_\ell^{-1}+u_m^{-1})
\end{multline}
\end{lemma}

This is equivalent to formula (4.5) of \cite{DFZJ07b}.
It is a slightly stronger version of the proposition of \cite{Zeil},
and can be proved along the same lines.
We use it to symmetrize the integrand:
\begin{multline}
S(L,p|t)
={1\over p!}\oint\cdots\oint
\prod_{\ell=1}^p
{\d u_\ell (1+t u_\ell)(1+\tau u_\ell)^{\tilde{p}+1} (\tau+u_\ell)^{\tilde{p}+1}\over 2\pi \i\,u_\ell}
\\
\prod_{1\leq \ell < m \leq p} (u_m-u_\ell)(\tau+u_\ell+u_m)(u_m^{-1}-u_\ell^{-1})(\tau+u_\ell^{-1}+u_m^{-1})
\end{multline}
Noting that $\prod_{1\leq \ell < m \leq p}
(u_m-u_\ell)(\tau+u_\ell+u_m)$ is just the Vandermonde determinant of
the $u_\ell(\tau+u_\ell)$ and similarly for the other factors, we can
finally pull the determinants out of the integral, resulting in:
\begin{equation}\label{Sint}
S(L,p|t)
=\det_{1\le\ell,m\le p}\left[
\oint {\d u \over 2\pi i u} (1+t u) u^{\ell - m+\tilde p} (\tau+u)^{\ell+\tilde p} (\tau+u^{-1})^{m+\tilde p}
\right]
\end{equation}
For general $t$, by using the binomial formula
we can evaluate this to be
\begin{equation}
S(L,p|t)
=\det_{1\le\ell,m\le p}\left[\sum_r
\tau^{\tilde{p}+2m+2\ell-2r-1}
{\ell+\tilde p\choose r-\ell}
\left(\tau{m+\tilde{p}\choose 2m-r}
+t {m+\tilde{p}\choose 2m-r-1}\right)
\right]
\end{equation}
where we recall that $\tilde{p}=L/2-1-p$.

At $t=\tau$ this expression simplifies slightly:
\begin{equation}\label{finSp}
S_+(L,p)=\det_{1\le\ell,m\le p}\left[\sum_r
\tau^{\tilde p+2m+2\ell-2r}
{\ell+\tilde p\choose r-\ell}
{m+\tilde p+1\choose 2m-r}
\right]
\end{equation}
as well as at $t=\tau^{-1}$:
\begin{equation}\label{finSm}
S_-(L,p)=\det_{1\le\ell,m\le p}\left[\sum_r
\tau^{\tilde p+2m+2\ell-2r}
{\ell+\tilde p+1\choose r-\ell}
{m+\tilde p\choose 2m-r}
\right]
\end{equation}
The summation over $r$ is such that only a finite number of terms is non-zero;
in practice a possible range is $0\le r\le 2p$.
One can check that formulae (\ref{finSp},\ref{finSm}) match
the expressions given in proposition~\ref{conj:S=T} for $L$ even.

\subsection{Odd size}
Assume $L=2n+1$. The reasoning being exactly identical to the case $L$ even,
we only provide the key formulae. The starting point is formally
the same integral formula as previously:
\begin{multline}\label{integinhomodd}
\psi_{a_1,\ldots,a_n}=
\prod_{1\leq i<j\leq L}
[1+x_i-x_j][1-x_i-x_j]
\oint\cdots\oint \prod_{\ell=1}^n
{\log q\over q-q^{-1}}{\d y_\ell\over \pi \i}\\
{\prod_{1\leq \ell<m\leq n}
[y_\ell-y_m][1+y_\ell-y_m][y_\ell+y_m]
\prod_{1\leq \ell\leq m\leq n} [1+y_\ell+y_m]
\over \prod_{\ell=1}^n\prod_{i=1}^{L} [y_\ell+x_i]
\prod_{i=1}^{a_\ell}[y_\ell-x_i]\prod_{i=a_\ell+1}^{L}
[1+y_\ell-x_i]
}
\end{multline}
but with an odd number of parameters $x_i$,
so that it produces a slightly different expression when the $x_i$ are set to zero:
\begin{multline}\label{integhomodd}
\psi_{a_1,\ldots,a_n}=\bar\psi_{b_1,\ldots,b_n}
= \oint\cdots\oint \prod_{\ell=1}^n {\d u_\ell\over 2\pi \i\, u_\ell^{b_\ell}} \left[
\prod_{1\leq \ell \leq m \leq n} (1-u_\ell u_m)(1+\tau u_m+u_\ell u_m)  \right.\\
\times \left. \prod_{1\leq \ell < m \leq n} (u_m-u_\ell)
(\tau+u_\ell+u_m) \right]
\end{multline}
Define, for $p$ a non-negative integer, $\tilde p=n-p$ and
\begin{equation}
S(L,p|t)=\sum_{\epsilon_1,\ldots,\epsilon_p\in \{0,1\}}
t^{\sum_{i=1}^p \epsilon_i}
\bar\psi_{1,\ldots,\tilde p,
\tilde p+2-\epsilon_1,\ldots,2n-\tilde p-\epsilon_p}
\end{equation}
so that $S_\pm(L,p)=S(L,p|\tau^{\pm 1})$.
We obtain the following integral expression for $S(L,p|t)$:
\begin{multline}
S(L,p|t)
= \oint\cdots\oint \left( \prod_{\ell=1}^{\tilde p} {\d u_\ell\over 2\pi \i\, u_\ell^\ell}\right)
\left( \prod_{\ell=\tilde{p}+1}^n {\d u_\ell (1+t u_\ell)\over 2\pi i u_\ell^{2\ell-\tilde{p}}}\right)
\left[ \prod_{1\leq \ell \leq m \leq n} (1-u_\ell u_m) \right.\\
\times \left. \prod_{\ell=1}^n (1+\tau u_\ell+u_\ell^2)
\prod_{1\leq \ell < m \leq n} (u_m-u_\ell)
(1+\tau u_m+u_\ell u_m) (\tau+u_\ell+u_m) \right]
\end{multline}
We integrate over the first $\tilde p$ variables and reindex the remaining ones:
\begin{multline}
S(L,p|t)
= \oint\cdots\oint
\prod_{\ell=1}^p {\d u_\ell (1+t u_\ell)(1+\tau u_\ell)^{\tilde{p}} (\tau+u_\ell)^{\tilde{p}}\over 2\pi \i\, u_\ell^{2\ell}}
\left[ \prod_{1\leq \ell \leq m \leq p} (1-u_\ell u_m) \right.\\
\times \left. \prod_{\ell=1}^n (1+\tau u_\ell+u_\ell^2)
\prod_{1\leq \ell < m \leq p} (u_m-u_\ell)
(1+\tau u_m+u_\ell u_m) (\tau+u_\ell+u_m) \right]
\end{multline}
We use lemma 2 and pull the determinants out of the integral as before:
\begin{align}
S(L,p|t)
&={1\over p!}\oint\cdots\oint
\prod_{\ell=1}^p
{\d u_\ell (1+t u_\ell)(1+\tau u_\ell)^{\tilde{p}} (\tau+u_\ell)^{\tilde{p}}
(1+\tau u_\ell+u_\ell^2)\over 2\pi \i\,u_\ell^2}
\notag\\
&\qquad\qquad\prod_{1\leq \ell < m \leq p} (u_m-u_\ell)(\tau+u_\ell+u_m)(u_m^{-1}-u_\ell^{-1})(\tau+u_\ell^{-1}+u_m^{-1})
\notag\\
&=\det_{1\le\ell,m\le p}\left[
\oint {\d u \over 2\pi \i\, u} (1+t u)(1+\tau u_\ell+u_\ell^2) u^{\ell - m+\tilde p-1} (\tau+u)^{\ell+\tilde p-1} (\tau+u^{-1})^{m+\tilde p-1}
\right]
\notag\\
&=\det_{1\le\ell,m\le p}\left[
\oint {\d u \over 2\pi \i\, u} (1+t u) u^{\ell - m+\tilde p} (\tau+u)^{\ell+\tilde p-1} (\tau+u^{-1})^{m+\tilde p}
\right]\label{Sintodd}
\end{align}
In the last line we wrote $(1+\tau u+u^2)/u=(\tau+u^{-1})+u$
and noted that the second term reproduces the column $m-1$ of the matrix
and thus can be subtracted without changing the determinant.

One finally obtains
\begin{equation}
S(L,p|t)
=\det_{1\le\ell,m\le p}\left[\sum_r
\tau^{\tilde{p}+2m+2\ell-2r-3}
{\ell+\tilde p-1\choose r-\ell}
\left(\tau{m+\tilde{p}\choose 2m-r}
+t {m+\tilde{p}\choose 2m-r-1}\right)
\right]
\end{equation}
where we recall that $\tilde p=(L-1)/2-p$.

At $t=\tau$ this simplifies to
\begin{equation}\label{finSpodd}
S_+(L,p)=\det_{1\le\ell,m\le p}\left[\sum_r
\tau^{\tilde p+2m+2\ell-2r-3}
{\ell+\tilde p-1\choose r-\ell}
{m+\tilde p+1\choose 2m-r}
\right]
\end{equation}
whereas at $t=\tau^{-1}$:
\begin{equation}\label{finSmodd}
S_-(L,p)=\det_{1\le\ell,m\le p}\left[\sum_r
\tau^{\tilde p+2m+2\ell-2r-1}
{\ell+\tilde p\choose r-\ell}
{m+\tilde p\choose 2m-r}
\right]
\end{equation}
thus reproducing the expressions
of proposition~\ref{conj:S=T} for $L$ odd.

\section{Proof of the bilinear recurrence relations for $T(L,p,k)$}
\label{HirotaProof}

In this section we prove Proposition \ref{prop-Trecur}. To simplify the formulae of this section we introduce $k'=L-2p-k$, so that definition \eqref{Tdef} reads
\be
T(L,p,k)\, =\, \det_{1\leq \ell, m\leq p}T_{\ell m}\, ,\qquad
T_{\ell m}(k,k') =\, \sum_r{\ell +k-1 \choose 2m-\ell-r}{m+k' \choose r}\tau^{2r}\, .
\label{Tdef-2}
\ee
Here the limits of summation in $r$ are automatically
fixed by the conditions ${n\choose r}=0~$  $\forall\,r<0$ and $\forall\,r>n$.
We begin with a derivation of yet another determinant formula for $T(L,p,k)$.

\begin{lemma}\label{T-U}
\be
T(L,p,k)\, =\, \det_{1\leq \ell, m\leq p+1}U_{\ell m}\, ,
\ee
where
\be
U_{\ell,1}(p) =\, (-1)^{\ell -1}\, \tau^{2(p+1-\ell)}\, , \qquad
U_{\ell, m+1}(k,k') =\, T_{\ell m}(k,k'-1)\, .
\label{Tdef-3}
\ee
\end{lemma}
\begin{proof}
To check the identity $\det U=\det T$ we shall perform linear transformations
of the matrix $U$ not affecting its determinant.
First, we combine adjacent rows of $U$ with the aim to set to zero all components of the first column, except the last element $U_{p+1,1}=(-1)^p$:
\[
V_{\ell m} = U_{\ell,{m+1}} + \tau^2 U_{\ell +1,m+1}.
\]
Then, by decomposing the determinant of the resulting matrix along the first column we find
\[
\det_{1\leq \ell, m\leq p+1}U_{\ell m}\,=\,\det_{1\leq \ell, m\leq p}V_{\ell m}\, .
\]

The lemma now follows by noticing that the rows of $T$ are linear combinations of those of $V$,
\begin{align*}
V_{\ell m}&=
\sum_r \left( {\ell +k-1 \choose 2m - \ell -r} {m+k'-1 \choose r} +
{\ell +k \choose 2m-\ell -r}{m+k'-1\choose r-1}\right)\tau^{2r}\\
&=
 \sum_r\left( {\ell +k-1\choose 2m - \ell -r}{m+k'\choose r}+
{\ell +k-1\choose 2m-\ell -r-1}{m+k'-1\choose r-1}\right)\tau^{2r}\\
&= T_{\ell m}(k,k')+\tau^2 T_{\ell,m-1}(k,k') ,
\end{align*}
where we have used Pascal's rule ${n\choose r} = {n-1\choose r}+ {n-1 \choose r-1}$ to go from the second to the third line. Hence we find
\[
\sum_{j=0}^{m-1} (-\tau^2)^j\, V_{\ell,m-j}\, =\, T_{\ell m},
\]
from which we conclude that $\det V =\det T$.
\end{proof}

Our derivation of formula \eqref{Trecur} is based on the use of the following particular example of
a Pl{\"u}cker relation for determinants (for the general case see
\cite{Sturmfels}). Consider a pair of $\, n\times n\, $ matrices $A _{\ell m}$
and $B_{\ell m}$. Denote by $A_\ell$ the $\ell$-th row of the matrix $A$
and  introduce notation
$$
\det A= |A|\, ,\qquad A\, =\,
\left[ A_1, \ldots, A_n \right],
$$
In this notation the Pl{\"u}cker relation  reads
\be
|A|\, |B| = \sum_{j=1}^{n}
\left| A_1 ,\ldots, A_{n-1},B_j \right|\times \left|
B_1, \ldots, B_{j-1},A_n,B_{j+1},\ldots, B_n \right|,
\label{Pluk}
\ee
where the sum is taken over permutations of the last row of $A$ with each row of $B$.

We now take $n=p+1$ and, recalling the definition of the matrix $U$ in \eqref{Tdef-3}, substitute for $A$ and $B$ the following matrices:
\begin{align}
A(k,k') &= U(k,k') = \left[ U_1 ,\ldots, U_{p+1} \right],
\label{A}\\
B(k,k') &= \left[ \delta_1, U_2 ,\ldots, U_p, \delta_{p+1} \right],
\label{B}
\end{align}
where  $\delta_{\ell m}$ is Kronecker's delta.

By  Lemma \ref{T-U} we have $|A|=T(L,p,k)$. To calculate the determinant of $B$ we expand along its top and bottom rows and then notice that the first column of the resulting $(p-1)\times (p-1)$ matrix is also of $\delta$-type:
$$
U_{\ell,2}(k,k') = T_{\ell,1}(k,k'-1) =\sum_r {\ell +k-1\choose 2-\ell-r}{k'\choose r} \tau^{2r}=
\delta_{\ell,2}\;\;\; \mbox{for}\;\; \ell\geq 2.
$$
So we calculate
\begin{align*}
|B|\, =\, \det_{1\leq\ell ,m\leq p-2} U_{\ell +2,m+2}(k,k')\,  &=\,
\det_{1\leq\ell ,m\leq p-2} T_{\ell +2,m+1}(k,k'-1)
\\
& =\, \det_{1\leq\ell ,m\leq p-2} T_{\ell m}(k+2,k')\, =\, T(L-2,p-2,k+2).
\end{align*}

Now let us consider the right hand side of the relation \eqref{Pluk}. Here only permutations
with the top and the bottom rows of the matrix $B$ give non-vanishing contributions.
Permuting the last row in $A$ and the first row in $B$ gives the following factors:
\begin{align*}
\left| U_1,\ldots U_p,\delta_1 \right| \,& =\, (-1)^p \det_{1\le \ell ,m\le p} U_{\ell,m+1}(k,k')
\\
&\hspace{5mm}=\,
(-1)^p \det_{1\le \ell ,m\le p} T_{\ell m}(k,k'-1)\, =\, (-1)^p\, T(L-1,p,k)\, ,
\end{align*}
\begin{align*}
&\left| U_{p+1},U_2,\ldots U_p, \delta_{p+1}\right| \, =\, (-1)^{p-1}
\left| U_2,\ldots,U_p,U_{p+1}\right|\hspace{21mm}
\\[1mm]
&\hspace{22mm} =\,
(-1)^p \det_{1\le \ell ,m\le p-1}U_{\ell m}(k+2,k'+1)\, =\,
(-1)^p\, T(L-1,p-2,k+2)\, ,
\end{align*}
where in the last calculation when passing to the second line
one i) expands the determinant along the first column noticing that
$U_{\ell +1,2} =\delta_{\ell,1}$ and ii) redefines matrix indices using the identity
$U_{\ell +2,m+1}(k,k')=U_{\ell m}(k+2,k'+1)\;\;\forall\;m>1$.

Permuting the last row in $A$ and the last row in $B$ gives the following factors:
\[
\left| U_1,\ldots, U_p,\delta_{p+1}\right| \,=\,
\left| U_1,\ldots,U_p\right|
=\,
\tau^2 \det_{1\le \ell ,m\le p} U_{\ell m}(k,k')\, =\,
\tau^2\, T(L-2,p-1,k)\, .
\]
The factor $\tau^2$ is extracted from the $p$-dependent first column of the matrix $U$:
$U_{\ell,1}(p+1) =\tau^2 U_{\ell,1}(p)$;
\begin{align*}
&\left| \delta_1,U_2,\ldots U_{p+1} \right| \, =\,  \det_{1\le \ell ,m\le p} U_{\ell +1,m+1}(k,k')
\\[1mm]
&\hspace{15mm} =\,
\det_{1\le \ell ,m\le p-1} U_{\ell +2,m+2}(k,k')\, =\,
\det_{1\le \ell ,m\le p-1} T_{\ell m}(k+2,k')
 \, =\,  T(L,p-1,k+2)\, .
\end{align*}
Here when passing to the second line of the calculation we first expand
$\det U_{\ell +1,m+1}$ along the first column $U_{\ell +1,1}=\delta_{\ell,1}$ and then
redefines indices of the matrix $T$: $T_{\ell +2,m+1}(k,k'-1)=T_{\ell m}(k+2,k')$.

Thus, the Pl\"{u}cker relation
\eqref{Pluk} for the matrices $A$ and $B$ defined in \eqref{A} and \eqref{B}
produces the equality \eqref{Trecur}.

\section{Proof of $\tau^2$-enumeration of ASMs}\label{ProofofPP2ASM}
In this section, we prove Proposition~\ref{PP2ASM}. We start from the determinant formulae of
\cite{Kuperberg} for the enumeration of various symmetry classes of Alternating Sign Matrices,
and reduce them to our own determinant formulae for $S_\pm$. In what follows,
we keep Kuperberg's notations (even though they are non-standard),
to ease the comparison of formulae with \cite{Kuperberg}. In each case, one starts with
configurations of the six-vertex model, which for certain particular boundary conditions
are identified with Alternating Sign Matrices in various symmetry classes.
There are three distinct Boltzmann weights for the six-vertex model,
taking into account $\mathbb{Z}_2$ symmetry, and they
are parametrised as $w_{\rm a}=ax^{-1}y-a^{-1}xy^{-1}$, $w_{\rm b}=a xy^{-1}-a^{-1}x^{-1}y$,
$w_{\rm c}=a^2-a^{-2}$, $x$ and $y$ being row/column spectral parameter and $a$ a global parameter.
In the ASM language, a weight $w_{\rm c}$ is assigned to a $\pm 1$, and weights
$w_{\rm b}$ and $w_{\rm a}$ are assigned to zeroes.
In the end we must take the homogeneous limit where $w_{\rm a}=w_{\rm b}$ and
$w_{\rm c}/w_{\rm a}=\tau=-q-q^{-1}$:
this ensures that adding a $-1$ to an ASM, that is {\it two}\/ extra vertices
of type c (adding a $-1$ also increases the number of $+1$ by 1), produces a weight $\tau^2$.
This is achieved by setting all spectral parameters to 1 and $a=-q$.
Similar parameters, called $b$ and $c$, which are related to boundary weights will be used below.

To be self-contained, we will give the matrices of Kuperberg relevant to this paper. It is useful to first define the functions $\sigma$ and $\alpha$ by
\[
\sigma(x)=x-x^{-1},\qquad \alpha(x) = \sigma(a x)\sigma(a/x).
\]
The relevant matrices then are
\begin{align}
M_{\rm U}(n;\vec{x},\vec{y})_{ij} &= \frac{1}{\alpha(x_i/y_j)} - \frac{1}{\alpha(x_iy_j)},
\label{MU}\\
M_{\rm UU}(n;\vec{x},\vec{y})_{ij} &= \frac{\sigma(b/y_j)\sigma(c x_i)}{\alpha(a x_i/y_j)} - \frac{\sigma(b/y_j)\sigma(c/x_i)}{\alpha(a/x_iy_j)} - \frac{\sigma(b y_j)\sigma(c x_i)}{\alpha(a x_i y_j)}
\nonumber\\
&{} + \frac{\sigma(b y_j)\sigma(c/x_i)}{\alpha(a y_j/x_i)},
\label{MUU}\\
M_{\rm HT}^{+}(n;\vec{x},\vec{y})_{ij} &= \frac{1}{\sigma(a y_j/x_i)} + \frac{1}{\sigma(a x_i/y_j)}.
\label{MHT}
\end{align}

Since the equalities of Proposition~\ref{PP2ASM} are known to be true at $\tau=1$
(they are consequences of the various relations between enumerations of ASMs found in
\cite{Kuperberg}, as well the relations between ASMs and PPs discussed in \cite{DFZJ07}),
and since both sides are easily checked to be polynomials in $\tau$ of the same degree,
we can safely drop various
trivial factors in the calculation, keeping only the determinant itself as well as factors that become
singular in the homogeneous limit.

\subsection{First formula}
\label{Uformula}
In \cite{Kuperberg}, it is explained how VSASMs are a special case of UASMs; more precisely,
the partition function of VSASMs can be obtained from the more general one of UASMs by tuning
a certain boundary parameter. It is also noted there that the boundary parameter only enters
the formula for the partition function in prefactors, and not in the determinant itself.
In particular the enumeration of VSASMs and UASMs are essentially the same.
We shall therefore write directly the partition function for UASMs of size $2n$
without all the regular prefactors; using \eqref{MU} it takes the form:
\begin{multline*}
Z_{\rm UASM}\propto \frac{1}{\Delta(x^2)\Delta(y^2)\Delta^*(x^2)\Delta^*(y^2)} \times\\
\det_{1\le i,j\le n}\left(\frac{1}{(a^2x_i^2-y_j^2)(a^2y_j^2-x_i^2)(a^2-x_i^2y_j^2)(1-a^2x_i^2y_j^2)}\right)
\end{multline*}
where $\Delta$ stands for the Vandermonde determinant, e.g.\ $\Delta(x^2)=\prod_{i<j}(x_j^2-x_i^2)$;
and $\Delta^*(x^2)=\prod_{i<j}(1-x_i^2x_j^2)$. The $x_i$, $y_j$ are spectral parameters which will eventually be set to one.

We use the following integral representation:
\begin{multline}
\frac{1}{(a^2x^2-y^2)(a^2y^2-x^2)(a^2-x^2y^2)(1-a^2x^2y^2)}
=
\\
\frac{a^3(a^2-1)^2}{\left(a^2+1\right) x^2y^4 (x^2-a^2) (a^2 x^2-1)}
\oint \frac{\d u}{2\pi \i}
\frac{u(u+a+a^{-1})}{\mu(u,x)\mu(u,1/x)\mu(u,a/y)\mu(u,ay)},
\label{intU}
\end{multline}
where
\be
\mu(u,x)=a(a+u)-(1+au)x^2.
\label{mudef}
\ee
The contour of integration in \eqref{intU} surrounds the $x$-dependent poles but not the $y$-dependent ones. This identity can be checked directly by residues. The various prefactors, as well as the integral sign, can be pulled out of the determinant and we thus find
\begin{multline*}
Z_{\rm UASM}\propto \frac{1}{\Delta(x^2)\Delta(y^2)\Delta^*(x^2)\Delta^*(y^2)}
\oint \frac{\d u_1(u_1+a+a^{-1})u_1}{2\pi \i}\cdots\oint \frac{\d u_n(u_n+a+a^{-1})u_n}{2\pi \i}\\
\det \left( \frac{1}{\mu(u_i,x_j) \mu(u_i,1/x_j)}\right)
\det \left( \frac{1}{\mu(u_i,a/y_j) \mu(u_i,ay_j)}\right)
\end{multline*}
In order to compute these determinants, we perform the following change of variables:
\[
X=-\frac{(1-x^2)^2}{(1-a^2x^2)(1-a^{-2}x^2)}
\qquad
Y=-\frac{(1-y^2)^2}{(1-a^2y^2)(1-a^{-2}y^2)}
\]
and use the factorizations
\begin{align*}
\mu(u,x) \mu(u,1/x)
&=(a^2-x^2)(a^2-x^{-2})(1-X u(a+a^{-1}+u)),
\\
\mu(u,ay) \mu(u,a/y)
&=a^2u^2(a^2-y^2)(a^2-y^{-2})(1-Y u^{-1}(a+a^{-1}+u^{-1})).
\end{align*}
Note that $a+a^{-1}=\tau$.
Again one can get rid of the trivial factors and obtain
\begin{multline*}
Z_{\rm UASM}\propto \frac{1}{\Delta(x^2)\Delta(y^2)\Delta^*(x^2)\Delta^*(y^2)}
\oint \frac{\d u_1(u_1+\tau)}{2\pi \i\, u_1}\cdots\oint \frac{\d u_n(u_n+\tau)}{2\pi \i\, u_n}\\
\det\left(\frac{1}{1-u_i(\tau+u_i)X_j}\right)
\det\left(\frac{1}{1-u_i^{-1}(\tau+u_i^{-1})Y_j}\right),
\end{multline*}
where the contours of integration surround the $X_i$-dependent poles.
The determinants are now of Cauchy type and can be evaluated exactly:
\begin{multline*}
Z_{\rm UASM}\propto \frac{\Delta(X)\Delta(Y)}{\Delta(x^2)\Delta(y^2)\Delta^*(x^2)\Delta^*(y^2)}
\oint \frac{\d u_1(u_1+\tau)}{2\pi \i\, u_1}\cdots\oint \frac{\d u_n(u_n+\tau)}{2\pi \i\, u_n}\\
\Delta(u(1+\tau u))\Delta(u^{-1}(1+\tau u^{-1}))
\prod_{i,j}\frac{1}{(1-u_i(\tau+u_i)X_j)(1-u_i^{-1}(\tau+u_i^{-1})Y_j)}.
\end{multline*}
The Vandermonde determinants outside the integral
cancel each other, leaving only a regular part, due to
\[
X_i-X_j=\frac{(1-a^2)^2(x_i^2-x_j^2)(1-x^2_i x_j^2)}{a^2(1-a^2x_i^2)(1-a^{-2}x_i^2)(1-a^2x_j^2)(1-a^{-2}x_j^2)}.
\]

At this stage one can take the homogeneous limit,
that is set $x_i=y_i=1$, or $X_i=Y_i=0$, which results in the simple expression
\[
A_{\rm V}(2n+1;\tau^2)\propto
\oint \frac{\d u_1(\tau+u_1)}{2\pi \i\, u_1}\cdots\oint \frac{\d u_n(\tau+u_n)}{2\pi \i\, u_n}
\Delta(u(\tau+ u))\Delta(u^{-1}(\tau+u^{-1})),
\]
where the integrals are performed around zero. We can once again exchange determinant and integral sign and find
\[
A_{\rm V}(2n+1;\tau^2)\propto
\det_{0\le\ell,m\le n-1} \oint \frac{\d u}{2\pi \i\, u} u^{\ell-m}(\tau+u)^{\ell+1}(\tau+u^{-1})^{m}
\]
The first column of this matrix being $(1,0,\ldots)$, we can restrict the range of the indices to
\[
A_{\rm V}(2n+1;\tau^2)\propto
\det_{1\le\ell,m\le n-1} \oint \frac{\d u}{2\pi \i\, u} u^{\ell-m}(\tau+u)^{\ell+1}(\tau+u^{-1})^{m}
\]
which is identical to \eqref{Sint} with $t=\tau^{-1}$, $p=n-1$, $\tilde p=0$.

\subsection{Second formula}
The reasoning is exactly the same; only the matrix elements of the determinant are slightly modified.
In the case of $A^{(2)}_{\rm VHPASM}$, we must use the matrix $M_{\rm UU}$ defined in \eqref{MUU} with
parameters $b=1/a$, $c=a$. In this case we use the following identity to represent the matrix elements in factorised form:
\begin{multline*}
\frac{-a^2 y^2(1+x^4)+\left(
\left(a^4+a^2+1\right)(1+y^4)-\left(a^2+1
   \right)^2 y^2
   \right) x^2
}{xy(a^2x^2-y^2)(a^2y^2-x^2)(a^2-x^2y^2)(1-a^2x^2y^2)}a^2
\\=(1-a^2)^2
\oint \frac{\d u}{2\pi \i}
\frac{u(1+u(a+a^{-1}))}{\mu(u,x)(u,1/x)
\mu(u,a/y)\mu(u,ay)}.
\end{multline*}

We thus find the following expression for the ``partition function'' $Z^{(2)}_{\rm VHPASM}$ (which is really a ratio of the partition function of UUASMs by the partition function of UASMs)
\begin{multline*}
Z^{(2)}_{\rm VHPASM}\propto \frac{\Delta(X)\Delta(Y)}{\Delta(x^2)\Delta(y^2)\Delta^*(x^2)\Delta^*(y^2)}
\oint \frac{\d u_1(1+\tau u_1)}{2\pi \i\,u_1}\cdots\oint \frac{\d u_n(1+\tau u_n)}{2\pi \i\,u_n}\\
\Delta(u(1+\tau u))\Delta(u^{-1}(1+\tau u^{-1}))
\prod_{i,j}\frac{1}{(1-u_i(\tau+u_i)X_j)(1-u_i^{-1}(\tau+u_i^{-1})Y_j)},
\end{multline*}
and in the homogeneous limit,
\begin{align*}
A^{(2)}_{\rm VHPASM}(2n+1;\tau^2)&\propto
\oint \frac{\d u_1(1+\tau u_1)}{2\pi \i\,u_1}\cdots\oint \frac{\d u_n(1+\tau u_n)}{2\pi \i\,u_n}
\Delta(u(\tau+ u))\Delta(u^{-1}(\tau+u^{-1}))\\
&\propto
\det_{0\le\ell,m\le n-1} \oint \frac{\d u}{2\pi \i} u^{\ell-m}(\tau+u)^{\ell}(\tau+u^{-1})^{m+1}.
\end{align*}
Once again one can remove the first line and column, and we recover
\eqref{Sint} with $t=\tau$, $p=n-1$, $\tilde p=0$.

\subsection{Third formula}
The generating function $\tilde{A}_{\rm UU}^{(2)}(n;\tau^2)$ is defined as one of the factors of $A_{\rm HT}^{(2)}(2n;\tau^2,1)$, the other being $A_{\rm UU}^{(2)}(n;\tau^2,1,1)$. To compute this generating function we need $Z_{\rm HT}^{+,(2)}(2n;(\vec{x},\vec{x}^{-1}),(\vec{y},\vec{y}^{-1}))$ which is defined in terms of the matrix $M_{\rm HT}^+$, see \eqref{MHT}. The matrix $M_{\rm HT}^+\left(2n;(\vec{x},\vec{x}^{-1}),(\vec{y},\vec{y}^{-1})\right)$ commutes with
\[
P = \begin{pmatrix} 0 & I_n \\ I_n & 0 \end{pmatrix},
\]
and thus can be brought to block-diagonal form with blocks
\[
M_{\rm HT}^{+,\pm}(n;\vec{x},\vec{y})_{ij} = M_{\rm HT}^+(n;\vec{x},\vec{y})_{ij} \pm M^+_{\rm HT}(n;\vec{x},\vec{y}^{-1})_{ij}.
\]
These matrices are readily identified with the matrix $M_{\rm UU}$ with $b=c=\i$ and $b=c=1$ respectively:
\begin{align*}
M_{\rm HT}^{+,-}(n;\vec{x},\vec{y})_{ij} &= \frac{1}{\sigma(\i x) \sigma(\i y)} M_{\rm UU}(n;\vec{x},\vec{y})_{ij}\qquad (b=c=\i),\\
M_{\rm HT}^{+,+}(n;\vec{x},\vec{y})_{ij} &= -\frac{1}{\sigma(x) \sigma(y)} M_{\rm UU}(n;\vec{x},\vec{y})_{ij}\qquad (b=c=1).
\end{align*}

The further reasoning is again analogous to that in Section~\ref{Uformula}, and in order to compute the homogeneous limit of $\det M_{\rm HT}^{(2)}(n;\vec{x},\vec{y})$ we now use the identity
\begin{multline}
\frac{
   (a^2+1)^2 x^2 y^2-a^2(x^2+y^2)(1+x^2y^2)
}{xy(a^2x^2-y^2)(a^2y^2-x^2)(a^2-x^2y^2)(1-a^2x^2y^2)}a^2
\\=(1-a^2)^2
\oint \frac{\d u}{2\pi\i}
\frac{u}{\mu(u,x)(u,1/x) \mu(u,a/y)\mu(u,ay)}.
\end{multline}
Thus we find
\begin{multline*}
\tilde Z^{(2)}_{\rm UU}\propto \frac{\Delta(X)\Delta(Y)}{\Delta(x^2)\Delta(y^2)\Delta^*(x^2)\Delta^*(y^2)}
\oint \frac{\d u_1}{2\pi \i\,u_1}\cdots\oint \frac{\d u_n}{2\pi \i\,u_n}\\
\Delta(u(1+\tau u))\Delta(u^{-1}(1+\tau u^{-1}))
\prod_{i,j}\frac{1}{(1-u_i(\tau+u_i)X_j)(1-u_i^{-1}(\tau+u_i^{-1})Y_j)},
\end{multline*}
and
\begin{align*}
\tilde A^{(2)}_{\rm UU}(4n;\tau^2)&\propto
\oint \frac{\d u_1}{2\pi \i\,u_1}\cdots\oint \frac{\d u_n}{2\pi \i\,u_n}
\Delta(u(\tau+ u))\Delta(u^{-1}(\tau+u^{-1}))\\
&\propto
\det_{0\le\ell,m\le n-1} \oint \frac{\d u}{2\pi \i} u^{\ell-m}(\tau+u)^{\ell}(\tau+u^{-1})^{m}.
\end{align*}
Removing the first line and column, we recover \eqref{Sintodd} with $t=\tau^{-1}$, $p=n-1$, $\tilde p=0$.


\newcommand\arxiv[1]{
\href{http://arxiv.org/abs/#1}{\tt arXiv:#1}}

\end{document}